\documentstyle{article}

\newcommand{\N}{{\rm I\kern-.5ex N}}
\newcommand{\Z}{{\sf \vrule height 1.55ex depth-1.2ex width.03em\kern-.11em Z
\kern-.9ex Z\kern-.11em\vrule height 0.3ex depth0ex width.03em}}
\newcommand{\Q}{{\rm\kern.2ex\vrule height1.55ex depth-.05ex width.03em\kern-.7ex Q}}
\newcommand{\R}{{\rm I\kern-.5ex R}}
\newcommand{\Rvar}{{\rm I\kern-.5ex R}}
\newcommand{\C}{{\rm\kern.3ex\vrule height1.55ex depth-.05ex width.03em\kern-.7ex C}}
\newcommand{\Cvar}{{\, \rm\kern.3ex\vrule height1.1ex depth-.05ex
width.03em\kern-.7ex C}}

\newcommand{\spat}{\hspace{4ex}}

\newcommand{\im}{\text{\rm Im}\,}
\newcommand{\re}{\text{\rm Re}\,}
\newcommand{\bld}{\text{\rm Ran}\,}
\newcommand{\lan}{\langle}
\newcommand{\ran}{\rangle}

\newcommand{\cB}{{\cal B}}
\newcommand{\cR}{{\cal R}}
\newcommand{\cL}{{\cal L}}

\newcommand{\cM}{{\cal M}}

\newcommand{\comp}{\,\rule[.5ex]{.2ex}{.2ex}\,}

\newcommand{\od}{\odot}
\newcommand{\ot}{\otimes}

\newcommand{\la}{\Lambda}

\newcommand{\om}{\omega}
\newcommand{\io}{\iota}

\newcommand{\vep}{\varepsilon}

\newcommand{\al}{\alpha}
\newcommand{\be}{\beta}

\newcommand{\gam}{\gamma}
\newcommand{\sde}{\delta}

\newcommand{\th}{\theta}

\newcommand{\oal}{\overline{\al}}

\newcommand{\text}[1]{\mbox{#1}}

\newcommand{\cst}{\text{C}$^*$}

\newcommand{\qed}{\ \hfill \rule{2mm}{2mm}}
\newenvironment{demo}{\medskip\noindent\bf Proof :\ \  \rm}{\qed\bigskip\par }

\newtheorem{definition}{Definition}[section]
\newtheorem{proposition}[definition]{Proposition}
\newtheorem{lemma}[definition]{Lemma}
\newtheorem{corollary}[definition]{Corollary}
\newtheorem{remark}[definition]{Remark}
\newtheorem{theorem}[definition]{Theorem}

\newtheorem{notation}[definition]{Notation}
\newtheorem{result}[definition]{Result}

\newtheorem{terminology}[definition]{Terminology}

\begin{document}
\begin{center}
\LARGE\bf
One-parameter representations on \cst-algebras
\end{center}

\bigskip

\begin{center}
\rm J. Kustermans  \footnote{Research Assistant of the
National Fund for Scientific Research (Belgium)}

Institut for Matematik og Datalogi

Odense Universitet

Campusvej 55

5230 Odense M

Denmark

\medskip

e-mail : johank@imada.ou.dk

\bigskip\medskip

\bf July 1997 \rm
\end{center}

\bigskip

\subsection*{Abstract}
Strongly continuous one-parameter representations on a \cst-algebra $A$ and their
extension to the multiplier algebra are investigated. We also give a proof of the Stone
theorem on Hilbert \cst-modules and look into some related problems.

\bigskip

\section*{Introduction}

In \cite{Zsido}, I. Cior\v{a}nescu and L. Zsid\'{o} investigated one-parameter groups on
Banach spaces with a certain well behaved weak topology on them. They looked mainly into
the correspondence between one-parameter groups (which are continuous with respect to this
weak topology) and their analytic generators.

In this paper, we look at the simpler case of one-parameter groups on a Banach space which
are continuous with respect to the norm topology. In most of the cases, we will even
restrict our attention to the case where this Banach space is a \cst-algebra or a
Hilbert-\cst-module.

\bigskip

The main result of this paper concerns one-parameter groups on \cst-algebras and their
extensions to the multiplier algebra :

\vspace{1mm}

Consider a \cst-algebra $A$ and a strongly continuous one-parameter group $\al$ on $A$
such that we have for every $t \in \R$ that $\al_t$ is a \cst-automorphism.

Then we have for every $t \in \R$ a unique \cst-automorphism $\oal_t$ on $M(A)$ which
extends $\al_t$. This gives us a strict strongly continuous one-parmeter group $\R
\rightarrow \text{Aut}(M(A)) : t \rightarrow \oal_t$.

Using the strict topology on $M(A)$, we can define also a strict analytic continuation of
this one-parameter group on $M(A)$. In this way, we get for every $z \in \C$ a linear
mapping $\oal_z$ in $M(A)$.

The strict closedness of this mapping $\oal_z$ turns out to be the only difficult result
to establish  (and this is not different from \cite{Zsido}, see theorem 2.4 of
\cite{Zsido}). The techniques used in \cite{Zsido} do not seem to work in this case (I do
not know whether the pair $(M(A),A^*)$ satisfies the axiom $\text{A}_1$ of \cite{Zsido}),
but the multiplicativity of each $\al_t$ allows us to circumvent this problem.

\vspace{1mm}

We will work in a setting which is a little bit more general. We only require that each
$\al_t$ satisfies some semi-multiplicativity property (and some useful one-parameter
groups satisfy this semi-mutltiplicativity property).

\bigskip

A motivation for studying such one-parameter groups can be found in the study of
\cst-algebraic quantum groups :
\begin{itemize}
\item It is probable that the left Haar weight is a \cst-algebraic KMS-weight.
So we have to work with the modular group of the left Haar weight which is  such a
\cst-algebraic one-parameter group.
\item At the moment, it seems also probable that the antipode of an algebraic
quantum group undergoes a polar decomposition : it can be written as the composition $R \,
\tau_{\frac{i}{2}}$ where $R$ is a $^*$-anti-automorphism of the \cst-algebra and $\tau$
is a \cst-algebraic one-parameter group. Therefore it is necessary to look at the
extension of $\tau_{\frac{i}{2}}$ to the multiplier algebra in order to investigate the
extension of the antipode.
\end{itemize}

\bigskip\medskip

In the first section, we will give an overview of the analytic continuation theory of
strongly continuous one parameter groups on Banach spaces. The second section looks into
one-parameter groups on \cst-algebras and their extensions to the multiplier algebra.

We quickly look into commuting one-parameter representations and tensor products of
one-parameter representations in section 3 and 4.

Section 5 revolves around a proof of the Stone theorem for Hilbert \cst-modules.

In section 6, we prove a familiar result concerning one-parameter representations
which are implemented by strictly positive elements in a Hilbert \cst-module.

The appendix contains an extension property for strict mappings between two \cst-algebras.

\bigskip

We end this introduction with some conventions and notations.

\medskip

The domain of any mapping $T$ will be denoted by $D(T)$, its range by $\text{Ran}\,T$.

\vspace{1mm}

All Banach spaces in this paper are considered over the complex numbers.
If $E$ is a Banach space, we denote the set of bounded operators by $\cB(E)$.

\vspace{1mm}

All Hilbert \cst-modules are considered as right modules and the inner products are linear
in the first variable. Standard references for Hilbert-\cst-modules are \cite{Baa1},
\cite{Baa}, \cite{Lan} and \cite{Wor6}. I have also collected some results in \cite{JK}.

Consider a Hilbert \cst-module $E$. Then we denote the set of adjointable operators on $E$
by $\cL(E)$. We denote the set of regular operators in $E$ by $\cR(E)$.

Let $T$ be a regular operator in $E$. Then we call $T$  positive $\Leftrightarrow$ $T$ is
self adjoint and $\lan T v, v \ran \geq 0$ for every $v \in D(T)$. We call $T$ strictly
positive $\Leftrightarrow$ $T$ is positive and has dense range.

If $T$ is strictly positive, then we can define for every $z \in \C$ the power $T^z$ which
is again a regular operator in $E$ (for more details, see \cite{JK}).

\bigskip

\section{One-parameter representations on Banach  spaces}

In this section, we will introduce the notion of one-parameter representations on a Banach
space and review the most important properties. Because it is a standard tool of
contemporary \cst-algebra theory, we will leave out some proofs but we will add some
comments. The standard reference for this section is \cite{Zsido}.

\bigskip

Whenever we use the notion of integrability of a function with values in a Banach space,
we mean the strong form of integrability (e.g. Analysis II, S. Lang) :

\medskip

Consider a measure space $(X,{\cal M},\mu)$ and a Banach space $E$.

\begin{itemize}
\item It is obvious how to define integrability for step functions from $X$ into $E$.
\item Let $f$ be a function from $X$ into $E$. Then $f$ is $\mu$-integrable if and
only if there exists a sequence of integrable step functions $(f_n)_{n=1}^\infty$ from $X$
into $E$ such that :
\begin{enumerate}
\item We have for almost every $x \in X$ that $(f_n(x))_{n=1}^\infty$ converges to $f(x)$
in the norm topology.
\item The sequence $(f_n)_{n=1}^\infty$ is convergent in the $L_1$-norm.
\end{enumerate}

In this case, the sequence $(\int f_n \, d\mu)_{n=1}^\infty$ is convergent and the
integral of $\int f \, d\mu$ is defined to be the limit of this sequence (Of course, one
has to prove that this limit is independent of the choice of the sequence
$(f_n)_{n=1}^\infty$).
\end{itemize}

\bigskip

We start this section with a basic but very useful lemma.

\begin{lemma} \label{o.lem1}
Consider Banach spaces $E$,$F$ and a closed linear mapping $\la$ from within $E$ into $F$.
Let $(X,\cM,\mu)$ be a measure space and $f$  a function from $X$ into $D(\la)$ such that
\begin{itemize}
\item $f$ is integrable.
\item The function $X \rightarrow F : t \mapsto \la(f(t))$ is integrable.
\end{itemize}
Then $\int f(t) \, d\mu(t)$ belongs to $D(\la)$ and $\la\bigl(\,\int f(t) \,
d\mu(t)\,\bigr)
= \int \la(f(t)) \, d\mu(t)$.
\end{lemma}
\begin{demo}
Define $G$ as the graph of the mapping $\la$. By assumption, we have that $G$ is a closed
subspace of $E \oplus F$. Next, we define the mapping $g$ from $X$ into $G$ such that
$g(t) = \bigl(f(t),\la(f(t))\bigr)$ for every $t \in X$. It follows that $g$ is integrable
and that
$$\int g(t) \, d\mu(t) = \bigl( \, \int f(t) \, d\mu(t) \, , \int \la(f(t)) \, d\mu(t)
\, \bigr) \ .$$
Because $G$ is a closed subspace of $E \oplus F$, we have that $\int g(t) \, d\mu(t)$
belongs to $G$. This implies that $\bigl(\, \int f(t) \, d\mu(t) \, , \int \la(f(t)) \,
d\mu(t) \,\bigr)$ belongs to $G$.
\end{demo}

\bigskip

When we speak of analyticity, we will always mean norm analyticity. But we have the
following well known result (which follows from the Uniform Boundedness principle, see 
e.g. the proof of theorem VI.4, Functional Analysis, Simon \& Reed).

\begin{result} \label{anal}
Consider a Banach space $E$ and a subspace $F$ of $E^*$ such that we have the equality
\newline $\|x\| = \sup\{ \, |\om(x)| \mid \om \in F \text{ with } \|\om\| \leq 1 \, \}$  \
for every $x \in E$. Let $O$ be an open subset of $\,\C$ and $f$ a function from $O$ into
$E$. Then $f$ is analytic $\Leftrightarrow$ We have for every $\om \in F$ that $\om \!
\circ \! f$ is analytic.
\end{result}

Notice that the result is true for $F=E^*$.

\bigskip\bigskip

Now we introduce the notion of strongly continuous one-parameter representations on a
Banach space.

\medskip

\begin{terminology}
Consider a Banach space $E$. By a one-parameter representation on $E$, we will always mean
a mapping $\al$ from $\R$ into $\cB(E)$ such that :
\begin{itemize}
\item We have for every $s,t \in \R$ that $\al_{s+t} = \al_s \al_t$.
\item $\al_0= \io$
\item We have for every $t \in \R$ that $\| \al_t \| \leq 1$.
\end{itemize}
We call $\al$  strongly continuous $\Leftrightarrow$ We have for every $a \in E$ that the
mapping $\R \rightarrow E : t \mapsto \al_t(a)$ is continuous.
\end{terminology}

This definition implies for every $t \in \R$ that $\al_t$ is invertible in $\cB(E)$, that
$(\al_t)^{-1} = \al_{-t}$ and that $\al_t$ is isometric.

\begin{remark}\rm We would like to mention the following special cases :
\begin{itemize}
\item If $E$ is a Hilbert \cst-module and $u$ is a strongly continuous one-parameter
representation on $E$ such that $u_t$ is a unitary element in $\cL(E)$ for every $t \in
\R$, we call $u$ a strongly continuous unitary one-parameter representation on $E$.
\item If $A$ is a \cst-algebra and $\al$ is a strongly continuous one-parameter
representation on $A$ such that $\al_t$ is a \cst-automorphism on $A$ for every $t \in
\R$, we call $\al$ a norm continuous one-parameter group on $A$.
\end{itemize}
\end{remark}

\bigskip\medskip

Consider $z \in \C$. Then $S(z)$ will denote the horizontal strip $ S(z) = \{ \, y \in \C
\mid \im y \in [0,\im z] \, \}$. We denote the interior of $S(z)$ by $S(z)^0$.

\bigskip

Now we will summarize the theory of analytic continuations of such strongly continuous
one-parameter representations. For the rest of this section, we will fix a Banach space
$E$ and a strongly continuous one-parameter representation $\al$ on $E$.

\begin{definition} \label{pa1.def1}
Consider $z \in \C$. We define the mapping $\al_z$ from within $E$ into $E$
such that :
\begin{itemize}
\item The domain of $\al_z$ is by definition the set \
$\{\, a \in E \mid$ There exists a function $f$ from $S(z)$ into $E$

such that
\begin{enumerate}
\item $f$ is continuous on $S(z)$
\item $f$ is analytic on $S(z)^0$
\item We have  that $\al_t(a) = f(t)$ for every $t \in \R$ \ \ \ \ \ \ $\}$
\end{enumerate}
\item Consider $a$ in the domain of $\al_z$ and $f$ the function from $S(z)$
into $E$ such that
\begin{enumerate}
\item $f$ is continuous on $S(z)$
\item $f$ is analytic on $S(z)^0$
\item We have  that $\al_t(a) = f(t)$ for every $t \in \R$
\end{enumerate}
Then we have by definition that $\al_z(a) = f(z)$.
\end{itemize}
\end{definition}

In the case that $z$ belongs to $\R$, this definition of $\al_z$ corresponds with $\al_z$
which we started from. Therefore the notation $\al_z$ is justified.

\begin{remark}\rm
\begin{itemize}
\item Consider $z \in \C$ and $y \in S(z)$. Then it is clear that $D(\al_z)$ is a
subset of $D(\al_y)$.
\item It is also clear that $D(\al_y)=D(\al_z)$ for $y,z \in \C$ with $\im  y =
\im z$.
\item Let $z \in \C$ and $a \in D(\al_z)$. Then the  function $S(z) \rightarrow E :
u \mapsto \al_u(a)$ is continuous on $S(z)$ and analytic on $S(z)^0$ (because this
function must be equal to the function $f$ from the definition).
\item Consider an element $a$ in $E$. We say that $a$ is analytic with respect to
$\al$ if $a$ belongs to $D(\al_z)$ for every $z \in \C$. If $a$ is analytic with respect
to $\al$, then the function $\C \rightarrow  E : u \mapsto \al_u(a)$ is analytic.
\end{itemize}
\end{remark}

\bigskip

We will now list the most important properties of the mappings $\al_z$. Most of the proofs
will be left out. Sometimes we will add some comments concerning the proofs. In most
cases, it comes down to a simple construction of an analytic function out of some
implicitly present analytic function(s). The techniques involved can be found in
\cite{Zsido}.

\bigskip

The proof of the following algebraic property is rather straightforward and only uses the
definition above.

\begin{proposition}
Consider $z \in \C$. Then the mapping $\al_z$ is a linear operator in $E$.
\end{proposition}

\medskip

\begin{result} \label{pa1.res1}
We have for every $z \in \C$, $t \in \R$ and $a \in D(\al_z)$ that  $\al_t(\al_z(a)) =
\al_{z+t}(a)$.
\end{result}

The proof of this result follows by noticing that the function $S(z) \rightarrow E : u
\mapsto \al_t(\al_u(a)) - \al_{u+t}(a)$ is continuous on $S(z)$, analytic on $S(z)^0$ and
equal to 0 on the real axis. Therefore it must be 0 on $S(z)$.

\medskip

\begin{result}
Consider $z \in \C$ and $a \in D(\al_z)$. Then the function $S(z) \rightarrow E : u
\mapsto \al_u(a)$ is bounded.
\end{result}

By continuity, the function above will be bounded on the set $[0, (\im z) \, i]$. The
previous result will therefore give us that our function is bounded on $S(z)$.

\medskip

Using the previous two results and the Phragmen-Lindelof theorem,
we get the following result.

\begin{corollary} \label{pa1.cor1}
We have for every $z \in \C$, $y \in S(z)$ and $a \in D(\al_z)$ that $\|\al_y(a)\|
\leq \max\{\|a\|,\|\al_z(a)\|\}$.
\end{corollary}

\medskip

In order to prove the next proposition, you prove, using only the definition, that $\al_y
\, \al_{-y}$ is the identity function on $D(\al_{-y})$ for $y \in \C$. Then the
proposition follows quickly.

\begin{proposition} \label{pa1.prop2}
Consider $z \in \C$. Then $\al_z$ is injective and $(\al_z)^{-1} = \al_{-z}$. Therefore
$\bld \al_z = D(\al_{-z})$ and $D(\al_z) = \bld \al_{-z}$.
\end{proposition}

\medskip

We already mentioned that $D(\al_{z+t})= D(\al_z)$ for every $z \in \C$ and $t \in \R$.
Combining this with the previous proposition, we arrive at the following conclusion.

\begin{result}
We have the equalities $D(\al_{z+t}) = D(\al_z)$ and $\bld \al_{z+t} =$ $\bld
\al_z$ for every $z \in \C$ and $t \in \R$.
\end{result}

\medskip

Using result \ref{pa1.res1}, we see that $\al_t \, \al_z = \al_{z+t}$ for every $z \in \C$
and $t \in \R$. Again using proposition \ref{pa1.prop2}, gives us the following
proposition.

\begin{proposition}
Consider $z \in \C$ and $t \in \R$. Then $\al_z \, \al_t = \al_t \, \al_z = \al_{z+t}$. So
$\al_t(D(\al_z)) = D(\al_z)$ and $\al_t(\bld \al_z) = \bld \al_z$.
\end{proposition}

\medskip

A generalization of the previous proposition is the following one. The proof of this
proposition consists of gluing two analytic functions together in the right way.

\begin{proposition} \label{pa1.prop4}
Consider $y,z \in \C$ such that $y$ and $z$ lie at the same side of the real axis. Then
$\al_y \, \al_z = \al_{y+z}$.
\end{proposition}

\medskip

Now we can use this proposition to prove a further generalization.

\begin{proposition}
Consider complex numbers $y$ and $z$. Then $\al_y \, \al_z \subseteq \al_{y+z}$. We have
moreover that $D(\al_y \, \al_z) = D(\al_z) \cap D(\al_{y+z})$ and that $\bld(\al_y \,
\al_z) = (\bld \al_y) \cap (\bld \al_{y+z})$.
\end{proposition}
\begin{demo}
We consider three cases.
\begin{enumerate}
\item $|\im z| \leq |\im y|$  and $y$ and $z$ lie on opposite sides of the
real axis :

Then $y+z$ and $-z$ lie at the same side of the real axis, so $\al_{y+z} \, \al_{-z} =
\al_y$  by the previous proposition. This implies that $\al_y \, \al_z \subseteq
\al_{y+z}$ and that $D(\al_y \, \al_z) = D(\al_z) \cap D(\al_{y+z})$.

\item $|\im z| \geq |\im y|$ and $y$ and $z$ lie on opposite sides of the real axis :

Then $y \in S(-z)$ which implies that $D(\al_{-z}) \subseteq D(\al_y)$. So we get that
$D(\al_y \, \al_z) = D(\al_z)$.

On the other hand, $y+z \in S(z)$. This implies that $D(\al_z) \subseteq D(\al_{y+z})$.
Consequently, $D(\al_z) \cap D(\al_{y+z}) = D(\al_z)$. So we see that $D(\al_y \, \al_z) =
 D(\al_z) = D(\al_z) \cap D(\al_{y+z})$.

We have furthermore that $y+z$ and $-y$ lie at the same side of the real axis, so
$\al_{-y} \, \al_{y+z} = \al_z$. Combining this with the result above, we get easily that
$\al_y \, \al_z \subseteq \al_{y+z}$.

\item $y$ and $z$ lie on the same side of the real axis :

Then we know already that $\al_{y+z} = \al_y \, \al_z$. Because $z \in S(y+z)$, we know
that $D(\al_{y+z}) \subseteq D(\al_z)$. Hence, $D(\al_z) \cap D(\al_{y+z}) = D(\al_{y+z})
= D(\al_y \, \al_z)$.
\end{enumerate}

Consequently we get in all the possible cases that $D(\al_y \, \al_z) = D(\al_z) \cap
D(\al_{y+z})$ and that $\al_y \, \al_z \subseteq \al_{y+z}$.

In a similar way, we get of course also that $D(\al_{-z} \, \al_{-y}) = D(\al_{-y})
\cap D(\al_{-(y+z)})$.
This implies that
$$(\bld \al_y) \cap (\bld \al_{y+z}) = D(\al_{-y}) \cap D(\al_{-(y+z)})
= D(\al_{-z} \, \al_{-y}) = \bld \bigl((\al_{-z} \, \al_{-y})^{-1}\bigr)
= \bld (\al_y \, \al_z) \ . $$
\end{demo}

\bigskip\medskip

As usual, we smear elements to construct elements which behave well with respect to $\al$.

\begin{notation}
Consider $a \in E$, $r > 0$ and $z \in \C$. Then we define the element $a(r,z)$ in $E$
such that
$$a(r,z) = \frac{r}{\sqrt{\pi}} \int \exp(-r^2 (t-z)^2) \, \al_t(a) \, dt
\ .$$
Then we have that $\|a(r,z)\| \leq \|a\| \, \exp(r^2 (\im z)^2)$.
\end{notation}

\medskip

For $a \in E$ and $r > 0$, we will use the notation
$$ a(r) = a(r,0) = \frac{r}{\sqrt{\pi}} \int \exp(-r^2 t^2) \, \al_t(a) \, dt \ .$$

\medskip

The following result is a direct consequence of the fact that the function
$$ \C \rightarrow E : u \mapsto \int \exp(-r^2 (t-u)^2) \, \al_t(a) \, dt $$
is analytic.

\begin{result} \label{pa1.prop1}
Let $a \in E$, $r > 0$ and $z \in \C$. Then $a(r,z)$ is analytic  with respect to $\al$
and $\al_y(a(r,z)) = a(r,z+y)$ for every $y \in \C$.
\end{result}

\medskip

\begin{corollary} \label{pa1.cor2}
We have for every $z \in \C$ that the mapping $\al_z$ is densely defined and has dense
range.
\end{corollary}

\medskip

The proof of the following result is an application of lemma \ref{o.lem1} and the fact
that $\al_t \, \al_y = \al_y \, \al_t$ for every $y \in \C$ and $t \in \R$.

\begin{result} \label{pa1.res2}
Consider   $r > 0$, $y,z \in \C$ and $a \in D(\al_y)$. Then $\al_y(a(r,z))$   $=
\al_y(a)(r,z)$.
\end{result}

\bigskip\bigskip

We will prove explicitly the following result to make a comparison with the more difficult
case of strict analytic continuations of norm continuous one-parameter groups on
\cst-algebras (see theorem \ref{pa2.th1} of the next section). We would also like to
mention that the proof if substantially easier than in the case of strongly continuous
one-parameter groups on von Neumann algebras (see \cite{Zsido}).

\begin{theorem} \label{pa1.th1}
Consider $z \in \C$. Then $\al_z$ is closed.
\end{theorem}
\begin{demo}
Choose a sequence $(a_k)_{k=1}^\infty$ in $D(\al_z)$, $a \in E$
and $b \in E$ such that $(a_k)_{k=1}^\infty$ converges to $a$
and $(\al_z(a_k))_{k=1}^\infty$ converges to $b$.

For every $k \in \N$, we define the function $f_k$ from $S(z)$ into $E$ such that $f_k(u)=
\al_u(a_k)$ for every $u \in S(z)$, then $f_k$ is bounded and continuous on $S(z)$ and
analytic on $S(z)^0$.

Choose $m,n \in \N$. We have for every $t \in \R$ that
\begin{itemize}
\item $\| f_m(t) - f_n(t) \|  = \| \al_t(a_m) - \al_t(a_n) \| \leq \|a_m - a_n \|$.
\item $\| f_m(z+t) - f_n(z+t) \| = \| \al_{z+t}(a_m) - \al_{z+t}(a_n) \|
$   $= \| \al_t(\al_z(a_m)) - \al_t(\al_z(a_n)) \| \leq \| \al_z(a_m) - \al_z(a_n)\| $.
\end{itemize}
Therefore, and because $f_m - f_n$ is bounded and continuous on $S(z)$ and analytic on
$S(z)^0$, the Phragmen Lindelof principle implies that
$$\| f_m(u) - f_n(u) \| \leq \max\{ \|a_m - a_n\| , \|\al_z(a_m)  - \al_z(a_n) \| \}
\text{\ \ \ \ \ \ \ \ (*)}$$
for every $u \in S(z)$.

This implies that the sequence $(f_k)_{k=1}^\infty$ is uniformly Cauchy on $S(z)$.
Consequently, there exists a function $f$ from $S(z)$ into $E$ such that
$(f_k)_{k=1}^\infty$ converges uniformly on $S(z)$. The uniform convergence implies that
$f$ is continuous on $S(z)$ and analytic on $S(z)^0$.

Letting $n$ tend to $\infty$ in inequality (*) gives that
$$\| f_k(u) - f(u) \| \leq \max\{ \|a_k - a\| , \|\al_z(a_k)  - b \| \}
$$ for every $u \in S(z)$ and $k \in \N$.
So $$\| \al_u(a_k) - f(u) \| \leq \max\{ \|a_k - a\| , \|\al_z(a_k)  - b \| \} $$
for every $u \in S(z)$ and $k \in \N$.

Choose $t$ in $\R$. The previous inequality implies that $(\al_t(a_k))_{k=1}^\infty$
converges to $f(t)$. On the other hand, it is also clear that $(\al_t(a_k))_{k=1}^\infty$
converges to $\al_t(a)$. This implies that $f(t)= \al_t(a)$

The previous inequality implies also that $(\al_z(a_k))_{k=1}^\infty$ converges to
$f(z)$.We also know that $(\al_z(a_k))_{k=1}^\infty$ converges to $b$. Hence, $f(z)$
equals $b$.

Consequently, we have by definition that $a$ belongs to $D(\al_z)$ and
$\al_z(a) = f(z) = b$.
\end{demo}

\bigskip

The proof of the following result is due to Woronowicz (see \cite{Wor3}). It guarantees the
existence of nice cores for $\al_z$.

\begin{result} \label{pa1.res3}
Consider a dense subset $K$ of $E$ and $R$ a subset of $\R^+_0$ such that $\infty$ is a
limit point of $R$. Define $C = \langle \, a(r) \mid a \in K, r \in R \, \rangle$ which
clearly consists of analytic elements with respect to $\al$.

Let $z$ be a complex number. Then $C$ is a core for $\al_z$.
\end{result}
\begin{demo}
Choose $x \in D(\al_z)$. Take $\vep > 0$.

Because $\infty$ is a limit point of $R$,  we get the existence of $r \in R$ such that
$$ \| x(r) - x \| \leq \frac{\vep}{2}   \hspace{1.5cm} \text{ and } \hspace{1.5cm}
\| (\al_z(x))(r) - \al_z(x) \| \leq \frac{\vep}{2} $$

Because $K$ is dense in $E$, there exists $a \in K$ such that $\|a - x\| \leq
\min \{ \frac{\vep}{2} ,  \frac{\vep}{2} \, \exp(- r^2 (\im z)^2) \}$.

This implies immediately that
$$
\|a(r) - x\|  \leq  \|a(r) - x(r)\| + \| x(r) - x \|
\leq  \| a - x \| + \frac{\vep}{2} \leq \frac{\vep}{2} + \frac{\vep}{2} = \vep
$$
Because $x$ belongs to $D(\al_z)$, we know that $\al_z(x(r)) = (\al_z(x))(r)$. So we get
also that
\begin{eqnarray*}
& & \|\al_z(a(r)) - \al_z(x)\| \leq  \|\al_z(a(r)) - \al_z(x(r))\| + \|\al_z(x(r)) -
\al_z(x) \|  \\
& & \spat =  \| a(r,z)  - x(r,z) \| + \|(\al_z(x))(r) - \al_z(x) \|  \leq
\exp( r^2 (\im z)^2 ) \, \|a - x \| + \frac{\vep}{2} \\
& & \spat \leq \exp( r^2 (\im z)^2 ) \, \exp(-r^2 (\im z)^2 ) \,
\frac{\vep}{2} + \frac{\vep}{2} = \vep
\end{eqnarray*}
\end{demo}

\medskip

\begin{corollary} \label{pa1.cor3}
Consider $z \in \C$ and $C$ a subspace of $D(\al_z)$ such that $\al_t(C) \subseteq C$ for
$t \in \R$ and such that $C$ is dense in $E$. Then $C$ is a core for $\al_z$.
\end{corollary}
\begin{demo}
It is clear that the mapping $C \rightarrow E : x \mapsto \al_z(x)$ is closable and we
define the closure of this mapping by $\be$. So $\be \subseteq \al_z$. Put $K = \lan \,
a(n) \mid a \in C, n \in \N \, \ran$.

\vspace{1mm}

Fix $n \in \N$ and $a \in C$ for a moment.

We have by assumption for every $t \in \R$ that $\al_t(a) \in C \subseteq D(\be)$ and that
$\be(\al_t(a)) = \al_z(\al_t(a)) = \al_t(\al_z(a))$. This implies that the function $\R
\rightarrow E : t \mapsto \exp(-n^2 t^2) \, \be(\al_t(a))$ is integrable. So lemma
\ref{o.lem1} implies that $a(n) \in D(\be)$.

\vspace{1mm}

So we see that $K \subseteq D(\be)$. The previous result implies that $K$ is a core for
$\al_z$, so we get that $\al_z = \be$. Hence $C$ is a core for $\al_z$.
\end{demo}

\bigskip

\begin{lemma} \label{pa1.lem1}
Consider $z \in \C$ and a core $C$ for $\al_z$. Let $x$ be an element in $D(\al_z)$. Then
there exists a sequence $(x_n)_{n=1}^\infty$ in $C$ such that
\begin{enumerate}
\item We have for every $n \in \N$ that $\|x_n\| \leq \|x\|$ and
$\|\al_z(x_n)\| \leq \|\al_z(x)\|$.
\item We have that $(x_n)_{n=1}^\infty$ converges to $x$ and that
$(\al_z(x_n))_{n=1}^\infty$ converges to $\al_z(x)$.
\end{enumerate}
\end{lemma}
\begin{demo} If $x=0$, the result above is trivially true.

So suppose that $x \neq 0$. Because $\al_z$ is injective, we have also that $\al_z(x) \neq
0$.

Because $C$ is a core for $\al_z$, there exists a sequence $(y_n)_{n=1}^\infty$ in $C$
such that $(y_n)_{n=1}^\infty$ converges to $x$ and such that $(\al_z(y_n))_{n=1}^\infty$
converges to $\al_z(x)$. Because $x \neq 0$ and $\al_z(x) \neq 0$, there exists $n_0 \in
\N$ such that $y_n \neq 0$ and $\al_z(y_n) \neq 0$ for every $n \in \N$ with $n \geq n_0$.

We define for every $n \in \N$ with $n \geq n_0$ the number $\lambda_n = \min \{
\frac{\|x\|}{\|y_n\|} , \frac{\|\al_z(x)\|}{\|\al_z(y_n)\|} \}$.

Then $(\lambda_n)_{n=n_0}^\infty$ converges to $1$.

We define for every $n \in \N$ with $n \geq n_0$ the element $x_n = \lambda_n \, y_n \in
C$. It is clear that $\|x_n\| \leq \|x\|$ and that $\|\al_z(x_n)\| \leq \|\al_z(x)\|$ for
every $n \in \N$ with $n \geq n_0$.

We have also that $(x_n)_{n=n_0}^\infty$ converges to $x$ and that
$(\al_z(x_n))_{n=n_0}^\infty$ converges to $\al_z(x)$.
\end{demo}

\bigskip

The next result guarantees that a weak continuity property on the continuation function
implies that it is automatically norm continuous.

\begin{proposition} \label{pa1.prop3}
Consider a subset $K$ of $E^*$ such that $K$ is separating for $E$ and $\th  \, \al_t \in
K$ for every $\th \in K$ and $t \in \R$. Let $z \in \C$ and $v,w \in E$ such that there
exists for every $\th \in K$ a function $f_\th$ from $S(z)$ in $\C$ such that
\begin{enumerate}
\item $f_\th$ is bounded and continuous on $S(z)$
\item $f_\th$ is analytic on $S(z)^0$
\item We have for every $t \in \R$ that $f_\th(t) = \th(\al_t(v))$
\item We have that $f_\th(z) = \th(w)$.
\end{enumerate}
Then $v$ belongs to $D(\al_z)$ and $\al_z(v) = w$.
\end{proposition}
\begin{demo} We have for every $n \in \N$ that $v(n) \in D(\al_z)$ and that
$$\al_z(v(n)) = \frac{n}{\sqrt{\pi}} \int \exp(-n^2 (t - z)^2) \, \al_t(v) \, dt$$

\medskip

Choose $\th \in K$ and $m \in \N$. Then we have that
$$\th\bigl(\al_z(v(m))\bigr)
= \frac{m}{\sqrt{\pi}} \int \exp(-m^2 (t - z)^2) \, \th(\al_t(v))  \, dt
= \frac{m}{\sqrt{\pi}} \int \exp(-m^2 (t - z)^2) \, f_\th(t)  \, dt $$

Because $f_\th$ is bounded and continuous on $S(z)$ and analytic on $S(z)^0$, the Cauchy
theorem implies that
$$\th\bigl(\al_z(v(m))\bigr) = \frac{m}{\sqrt{\pi}} \int \exp(-m^2 t^2) \, f_\th(t+z)
\, dt \hspace{1cm} (*)$$

\medskip

Take $s \in \R$. Then we have for every $u \in \R$ that
$$f_\th(u+s) - f_{\th \al_s}(u) = \th(\al_{u+s}(v)) - (\th \al_s)(\al_u(v))
= \th\bigl(\al_s(\al_u(v))\bigr) - \th\bigl(\al_s(\al_u(v))\bigr) = 0$$
Because the function $S(z)  \rightarrow \C : u \mapsto f_\th(u+s) - f_{\th \al_s}(u)$ is
continuous on $S(z)$ and analytic on $S(z)^0$, this implies that $$f_{\th}(z+s) = f_{\th
\al_s}(z) = (\th \al_s)(w) = \th(\al_s(w))$$

\medskip

Substituting this equality in (*), we get that
$$\th\bigl(\al_z(v(m))\bigr) = \frac{m}{\sqrt{\pi}} \int \exp(-m^2 t^2) \, \th(\al_t(w))
\, dt = \th(w(m))$$

Because $K$ is separating for $E$, we see that
$\al_z(v(n)) = w(n)$ for every $n \in \N$.

\medskip

From this all, we get that $(v(n))_{n=1}^\infty$ converges to $v$ and that
$\bigl(\al_z(v(n))\bigr)_{n=1}^\infty$ converges to $w$.
Hence the closedness of $\al_z$ implies that $v \in D(\al_z)$ and
$\al_z(v) = w$.
\end{demo}

\bigskip

\section{Extensions of analytic continuations of a one-parameter group}

In this section, we will fix a \cst-algebra $A$. We will look into strongly continuous
one-parameter representations which have a natural extension to strictly continuous
one-parameter representations on $M(A)$. We will prove that these strict extensions give
rise to strictly analytic continuations.

\bigskip

\begin{definition}
Consider a strongly continuous one-parameter representation $\al$ on $A$. We call $\al$
semi-multiplicative if and only if there exist strongly continuous one-parameter
representations $\be$ and $\gam$ on $A$ such that
$$\be_t(b) \, \al_t(a) = \al_t(b\,a) \hspace{1.5cm} \text{and} \hspace{1.5cm} \al_t(a)
\, \gam_t(b) = \al_t(a \, b)$$
for every $a,b \in A$ and $t \in \R$.
\end{definition}

\medskip

It is clear that $\be$ and $\gam$ in the definition above are uniquely determined.
Therefore we introduce the following notation.

\begin{notation}
Consider a strongly continuous one-parameter representation $\al$ on $A$ which is
semi-multiplicative. Then we define the strongly continuous one-parameter representations
$\al^l$ and $\al^r$ on $A$ such that
$$\al^l_t(b) \, \al_t(a) = \al_t(b\,a) \hspace{1.5cm} \text{and} \hspace{1.5cm} \al_t(a)
\, \al^r_t(b) = \al_t(a \, b)$$
for every $a,b \in A$ and $t \in \R$.
\end{notation}

\medskip

\begin{remark} \rm
Of course, a natural case arises when we have a strongly continuous one-parameter
representation $\al$ on $A$ such that $\al_t$ is multiplicative for every $t \in \R$. Then
$\al$ is semi-multiplicative and $\al^l_t = \al^r_t = \al_t$ for every $t \in \R$.
\end{remark}

\bigskip\medskip

For the rest of this section, we consider a strongly continuous one-parameter
representation $\al$ on $A$ which is semi-multiplicative.

\medskip

Consider $t \in \R$. Then it is not difficult to check that $\al^l_t$ and $\al^r_t$ are
multiplicative. This implies that $\al^l$ and $\al^r$ satisfy the same condition as $\al$.

\medskip

It is also easy to see that $\al_t$, $\al^l_t$ and $\al^r_t$ are strictly continuous. By
proposition \ref{app.prop1}, this gives us strict linear mappings $\overline{\al_t}$,
$\overline{\al^l_t}$ and $\overline{\al^r_t}$.

\medskip

We have moreover that
$$\al^l_t(b) \, \al_t(a) = \al_t(b \, a) \hspace{1.5cm} \text{and} \hspace{1.5cm}
\al_t(a) \, \al^r_t(b) = \al_t(a \, b)$$
for every $a,b \in M(A)$.

\medskip

This implies that $\overline{\al_t}$, $\overline{\al^l_t}$ and $\overline{\al^r_t}$ are
strictly continuous (and not merely strictly continuous on bounded sets).

\bigskip

By the results of the previous section , we have for any $z \in \C$ a closed operator
$\al_z$ from within $A$ into $A$ which is generally unbounded.

Using the strict topology on $M(A)$, we will construct an extension of $\al_z$ to a linear
operator from within $M(A)$ into $M(A)$. The most significant result in this section
(theorem \ref{pa2.th1}) says that this extension is strictly closed.

\bigskip

Concerning analyticity, we have the following result. The proof of this fact uses result
\ref{anal} and the fact that every continuous linear functional on $A$ is of the form $a
\, \om$ with $\om \in A^*$ and $a \in A$.

\begin{result}
Consider an open subset $O$ of the complex plane, $f$ a function from $O$ into $M(A)$.
\newline  Then $f$ is analytic $\Leftrightarrow$ We have for every $\om \in A^*$ that the
function $\overline{\om} \! \circ \! f$ is analytic $\Leftrightarrow$ We have for every $a
\in A$ that the function $O \rightarrow A : z \mapsto f(z) \, a$ is analytic.
\end{result}

\bigskip

The proof of the following algebraic result is not difficult and will therefore be left
out. It involves only the construction of an analytic function out of the obvious
implicitly present ones.

\begin{proposition}
Consider a complex number $z$ and $a \in D(\al_z)$. Then we have the following properties
:
\begin{itemize}
\item We have for every $b \in D(\al^l_z)$ that $b\,a$ belongs to $D(\al_z)$ and
$\al_z(b \, a) = \al^l_z(b) \al_z(a)$.
\item We have for every $b \in D(\al^r_z)$ that $a\,b$ belongs to $D(\al_z)$ and
$\al_z(a \, b) = \al_z(a) \al^r_z(b)$.
\end{itemize}
\end{proposition}

\bigskip

For the rest of this section, we want to concentrate on strictly analytic continuations of
$\al$. All the results concerning norm analytic continuations will remain true (or have an
obvious strict variant). Most of the proofs can be easily translated to this case and will
therefore be omitted. A major exception to this rule will be the proof of the strict
closedness of the strict analytic continuation. The proof of this result will be a little
bit more involved than the proof of theorem \ref{pa1.th1}.

\bigskip

The next result is easy to check.

\begin{result}
We have for every $s,t \in \R$ and $a \in M(A)$ that
$\al_s(\al_t(a)) = \al_{s+t}(a)$.
\end{result}

\begin{proposition} \label{pa2.prop1}
Let $a$ be an element in $M(A)$. Then the function $\R \rightarrow
M(A) :  t \mapsto \al_t(a)$ is strictly continuous.
\end{proposition}
\begin{demo}
Choose $t \in \R$ and a sequence $(t_n)_{n=1}^\infty$ in $\R$ such that
$(t_n)_{n=1}^\infty$ converges to $t$.

Choose $b \in A$. Take $\vep >0$. Then we have that
\begin{itemize}
\item There exists $n_1 \in \N$ such that $\|\al^r_{t_n - t}(b) - b\|
\leq \frac{\vep}{2(\|a\|+1)}$ for every $n \in \N$ with $n \geq n_1$.
\item There exists $n_2 \in \N$ such that $\|\al_{t_n - t}( \, \al_t(a) b \,)
- \al_t(a) b \| \leq \frac{\vep}{2}$ for every $n \in \N$ with $n \geq n_2$.
\end{itemize}

Put $n_0 = \max\{n_1,n_2\}$. Then we have for every $n \in \N$ with $n \geq n_0$ that
\begin{eqnarray*}
& & \|\al_{t_n}(a) b - \al_t(a) b \| =  \|\al_{t_n - t}(\al_t(a)) b - \al_t(a) b \| \\ & &
\spat \leq \|\al_{t_n - t}(\al_t(a)) b - \al_{t_n - t}(\al_t(a)) \al^r_{t_n-t}(b)\| + \|
\al_{t_n - t}(\al_t(a)) \al^r_{t_n - t}(b) - \al_t(a) b \| \\ & & \spat \leq \|\al_{t_n -
t}(\al_t(a))\| \, \|b - \al^r_{t_n -t }(b)\| + \|\al_{t_n - t}(\al_t(a) b)  -
\al_t(a) b \| \\ & & \spat \leq \|a\| \, \frac{\vep}{2(\|a\|+1)} + \frac{\vep}{2} \leq
\vep \ .
\end{eqnarray*}

So we see that $(\al_{t_n}(a) b)_{n=1}^\infty$ converges to $\al_t(a) b$.

\medskip

Using $\al^l$ in stead of $\al^r$, we get  in a similar way that $(b \, \al_{t_n}(a)
)_{n=1}^\infty$ converges to $b \, \al_t(a)$.
\end{demo}

An obvious generalization of definition \ref{pa1.def1} is the following one.

\begin{definition} \label{pa2.def1}
Consider $z \in \C$. We define the mapping $\oal_z$ from within $M(A)$ into $M(A)$ as
follows.
\begin{itemize}
\item The domain of $\oal_z$ is by definition the set \
$\{\, a \in M(A) \mid$ There exists a function $f$ from $S(z)$ into

$M(A)$ such that
\begin{enumerate}
\item $f$ is strictly continuous on $S(z)$
\item $f$ is analytic on $S(z)^0$
\item We have  that $\al_t(a) = f(t)$ for every $t \in \R$ \ \ \ \ \ \ $\}$
\end{enumerate}
\item Choose $a$ in  the domain of $\oal_z$ and let $f$ be the function from $S(z)$
into $M(A)$ such that
\begin{enumerate}
\item $f$ is strictly continuous on $S(z)$
\item $f$ is analytic on $S(z)^0$
\item We have  that $\al_t(a) = f(t)$ for every $t \in \R$
\end{enumerate}
Then we have by definition that $\oal_z(a) = f(z)$.
\end{itemize}
\end{definition}

\medskip

In the case where $z$ belongs to $\R$, the previous proposition implies that  $\oal_z$ is
equal to $\overline{\al_z}$ (as defined in proposition \ref{app.prop1}). It is also clear
that $\al_z \subseteq \oal_z$ for every $z \in \C$. These two previous remarks justify the
following notation.

\begin{notation}
For every $z \in \C$ and $a \in D(\oal_z)$, we put $\al_z(a) = \oal_z(a)$.
\end{notation}

\medskip

\begin{remark}\rm
\begin{itemize}
\item Consider $z \in \C$ and $y \in S(z)$. Then it is clear that $D(\oal_z)$ is a subset
of $D(\oal_y)$.
\item It is also clear that $D(\oal_y)=D(\oal_z)$ for $y,z \in \C$ with $\im y =
\im z$.
\item Let $z \in \C$ and $a \in D(\oal_z)$. Then the  function $S(z) \rightarrow M(A) :
u \mapsto \al_u(a)$ is strictly continuous on $S(z)$ and analytic on $S(z)^0$ (because
this function must be equal to the function $f$ from the definition).
\item Consider an element $a$ in $M(A)$. We say that $a$ is strictly analytic with respect
to $\al$ if $a$ belongs to $D(\oal_z)$ for every $z \in \C$. If $a$ is strictly analytic
with respect to $\al$, then the function $\C \rightarrow  M(A) : u \mapsto \al_u(a)$ is
analytic.
\end{itemize}
\end{remark}

Again, the proof of the linearity is straightforward, just by using
the definition above.

\begin{proposition}
The mapping $\oal_z$ is a linear operator in $M(A)$.
\end{proposition}

\medskip

\begin{result} \label{pa2.res1}
We have fore every $z \in \C$, $t \in \R$ and $a \in D(\oal_z)$ that $\al_t(\al_z(a)) =
\al_{z+t}(a)$.
\end{result}
\begin{demo}
Define the function $f$ from $S(z)$ into $M(A)$ such that $f(u) = \al_{u+t}(a) -
\al_t(\al_u(a))$ for every $u \in S(z)$. Because $\overline{\al_t}$ is strictly continuous,
we get that $f$ is strictly continuous on $S(z)$. It is also clear that $f$ is analytic on
$S(z)^0$ and that $f=0$ on $\R$. As usual, this implies that $f=0$ on $S(z)$. In
particular, $f(z)=0$.
\end{demo}

\begin{result}
Consider $z \in \C$ and $a \in D(\al_z)$. Then the function $S(z) \rightarrow M(A) : u
\mapsto \al_u(a)$ is bounded.
\end{result}
\begin{demo}
Choose $b \in A$. Then the function $[0, (\im z) i ] \rightarrow A : u \mapsto
\al_u(a) b$ is continuous. Therefore, it is bounded. So there exists a positive number
$M_b$ such that $\| \al_u(a) b\| \leq M_b$ for every $u  \in [0, (\im z) i ]$.

Using the uniform boundedness principle, we get the existence of a positive number $M$
such that $\|\al_u(a)\| \leq M$ for every $u \in [0, (\im z) i ]$. This implies for every
$u \in S(z)$ that
$$\|\al_u(a)\| = \|\al_{\text{\scriptsize Re}\,u + i\,\text{\scriptsize Im}\,u}(a) \|
\stackrel{(*)}{=} \| \al_{\text{\scriptsize Re}\,u}(\, \al_{i\,\text{\scriptsize Im}\,u}(a)
\,) \| = \| \al_{i\,\text{\scriptsize Im}\,u}(a)\| \leq M  \  ,$$
where we used the previous result in equation (*).
\end{demo}

\medskip

\begin{corollary}
Consider $z \in \C$ and $a \in D(\oal_z)$. Then $\|\al_y(a)\| \leq
\max\{\|a\|,\|\al_z(a)\|\}$ for every $y \in S(z)$.
\end{corollary}
\begin{demo}
Choose $y \in S(z)$. Take $\om \in A^*$ with $\|\om\| \leq 1$.

The function $S(z) \rightarrow \C : u \mapsto \om(\al_u(a))$ is continuous and bounded on
$S(z)$ and analytic on $S(z)^0$.

We have for every $t \in \R$ that
$$|\om(\al_t(a))| \leq  \|\al_t(a)\| \leq \|a\| \hspace{1cm} \text{and} \hspace{1cm}
|\om(\al_{z+t}(a))| \leq \|\al_{z+t}(a)\| =  \|\al_t(\al_z(a))\|  \leq \|\al_z(a)\| $$
Hence, Phragmen Lindelof implies that $|\om(\al_y(a))| \leq \max\{\|a\|,\|\al_z(a)\|\}$.

So we get that $\|\al_y(a)\| \leq \max\{\|a\|,\|\al_z(a)\|\}$.
\end{demo}

\bigskip

The next proposition is an easy consequence of proposition \ref{pa1.prop3}.

\begin{proposition} \label{pa2.prop9}
Consider $z \in \C$ and $a \in A \cap D(\oal_z)$. Then $a$ belongs to
$D(\al_z)$ $\Leftrightarrow$ $\al_z(a)$ belongs to $A$.
\end{proposition}

The proof of the next proposition is rather straightforward and will therefore be left
out.

\begin{proposition} \label{pa2.prop3}
Let $z$ be a complex number and $a \in D(\oal_z)$. Then we have the following properties :
\begin{itemize}
\item We have for every $b \in D(\al^l_z)$ that $b a$ belongs to $D(\al_z)$ and
$\al_z(b a) = \al^l_z(b) \al_z(a)$.
\item We have for every $b \in D(\al^r_z)$ that $a b$ belongs to $D(\al_z)$ and
$\al_z(a b) = \al_z(a) \al^r_z(b)$.
\end{itemize}
\end{proposition}

\bigskip

The proof of the next results are just the same as in the case of norm analytic
coninuations.

\begin{proposition} \label{pa2.prop7}
Consider $z \in \C$. Then $\oal_z$ is injective and $(\oal_z)^{-1} =   \oal_{-z}$.
Therefore $\bld \oal_z = D(\oal_{-z})$ and $D(\oal_z) = \bld
\oal_{-z}$.
\end{proposition}

\medskip

We already mentioned that $D(\oal_{z+t})= D(\oal_z)$ for every $z \in \C$ and $t \in \R$.
Combining this with the previous proposition, we arrive at the following conclusion.

\begin{result}
We have the equalities $D(\oal_{z+t}) = D(\oal_z)$ and $\bld \oal_{z+t} =$ $\bld \oal_z$
for every $z \in \C$ and $t \in \R$.
\end{result}

\medskip

Using result \ref{pa2.res1}, we see that $\oal_t \, \oal_z = \oal_{z+t}$ for every $z \in
\C$ and $t \in \R$. Using proposition \ref{pa2.prop7}, this gives us the following
proposition.

\begin{proposition}
Consider $z \in \C$ and $t \in \R$. Then $\oal_z \, \oal_t = \oal_t \, \oal_z =
\oal_{z+t}$. So $\oal_t(D(\oal_z)) = D(\oal_z)$ and $\oal_t(\bld \oal_z) = \bld
\oal_z$.
\end{proposition}

\medskip

A generalization of the previous proposition is the following one.

\begin{proposition}
Consider $y,z \in \C$ such that $y$ and $z$ lie at the same side of the real axis. Then
$\oal_y \, \oal_z = \oal_{y+z}$.
\end{proposition}

\medskip

As for norm analytic continuations, this is then used to prove the next generalization :

\begin{proposition}
Consider complex numbers $y$ and $z$. Then $\oal_y \, \oal_z \subseteq
\oal_{y+z}$. We have moreover that $D(\oal_y \, \oal_z) = D(\oal_z) \cap D(\oal_{y+z})$
and $\bld(\oal_y \, \oal_z) = (\bld \oal_y) \cap (\bld \oal_{y+z})$.
\end{proposition}

\medskip

Also in this case can we smear elements. This is possible because of proposition
\ref{pa2.prop1}.

\begin{notation} \label{pa2.not1}
Consider $a \in M(A)$, $r > 0$ and $z \in \C$. Then we define the element $a(r,z)$ in
$M(A)$ such that
$$a(r,z) \, b = \frac{r}{\sqrt{\pi}} \int \exp(-r^2 (t-z)^2) \, \al_t(a)b \, dt $$
for every $b \in A$. Then we have that $\|a(r,z)\| \leq \|a\| \, \exp(r^2 (\im z)^2)$.
\end{notation}

\medskip

For $a \in M(A)$ and  $r > 0$, we will use the notation $a(r) = a(r,0)$. So $$a(r) \, b =
\frac{r}{\sqrt{\pi}} \int \exp(-r^2 t^2) \, \al_t(a) b \, dt $$ for every $b \in A$.

\medskip

\begin{result}
Let $a \in M(A)$, $r > 0$ and $z \in \C$. Then $a(r,z)$ is strictly analytic with respect
to $\al$ and $\al_y(a(r,z)) = a(r,z+y)$ for every $y \in \C$.
\end{result}

\medskip

\begin{result}
Consider   $r > 0$, $y,z \in \C$ and $a \in D(\oal_y)$. Then $\al_y(a(r,z)) =
\al_y(a)(r,z)$.
\end{result}

We should be careful by just saying that we can copy the proof of result \ref{pa1.res2} to
prove this result. Notice that in the proof of \ref{pa1.res2}, we use lemma \ref{o.lem1}
which concerns the norm and not the strict topology. But this can be easily circumvented.

Choose $b \in D(\al^r_y)$. Then lemma \ref{o.lem1} and the fact that $\oal_y \, \oal_t =
\oal_t \, \oal_y$ for every $t \in \R$ imply that $a(r,z) \, b$ belongs to $D(\al_y)$
and $\al_y(a(r,z) \, b) = \al_y(a)(r,z) \, \al^r_y(b)$. So we get that $\al_y(a(r,z)) \,
\al^r_y(b) = \al_y(a)(r,z) \, \al^r_y(b)$.

This implies that $\al_y(a(r,z)) = \al_y(a)(r,z)$.

\bigskip

Another useful property of smearing is contained in the following result. The proof of
this result is not difficult but involves some fiddling around. A similar result can be
found in \cite{Val} (proposition 0.3.2). The proof there uses Ascoli's theorem, whilst the
proof here depends on basic integration techniques.

\begin{proposition} \label{pa2.prop6}
Let $z \in \C$ and $r > 0$. Consider $a \in M(A)$ and a bounded net $(a_i)_{i \in I}$ in
$M(A)$ such that $(a_i)_{i \in I}$ converges strictly to a. Then $(a_i(r,z))_{i \in I}$ is
also bounded and converges strictly to $a(r,z)$.
\end{proposition}
\begin{demo}
By notation \ref{pa2.not1}, we see immediately that $(a_i(r,z))_{i \in I}$ is bounded.

Choose $b \in A$. Put $c = \al^r_{-\text{\scriptsize Re}\,z}(b)$.

We have for every $i \in I$ that
\begin{eqnarray*}
\|a_i(r,z)\,b - a(r,z)\,b \|
& = & \| \frac{r}{\sqrt{\pi}} \int \exp(-r^2 (t-z)^2) \, \al_t(a_i)b \, dt
- \frac{r}{\sqrt{\pi}} \int \exp(-r^2 (t-z)^2) \, \al_t(a) b  \, dt  \| \\
& \leq & \frac{r}{\sqrt{\pi}} \int |\exp(-r^2(t-z)^2)| \,
\|\al_t(a_i) b - \al_t(a) b \| \, dt \\
& = & \frac{r}{\sqrt{\pi}} \exp(r^2 (\im z)^2) \int \exp(-r^2(t -\re z)^2)
\,  \|\al_t(a_i \, \al^r_{-t}(b) - a \, \al^r_{-t}(b)) \| \, dt  \\
& = & \frac{r}{\sqrt{\pi}} \exp(r^2 (\im z)^2) \int \exp(-r^2(t -\re z)^2)
\,  \|a_i \, \al^r_{-t}(b) - a \, \al^r_{-t}(b) \| \, dt  \\
& = & \frac{r}{\sqrt{\pi}} \exp(r^2 (\im z)^2) \int \exp(-r^2 t^2)  \,
\|(a_i-a) \al^r_{-t}(c) \| \, dt
\end{eqnarray*}
So we see that in order to prove that $(a_i(r,z)\,b)_{i \in I}$ converges to $a(r,z)\,b$,
it is sufficient to prove that the net
$$\left( \, \,  \int \exp(-r^2 t^2) \, \|(a_i-a)\al^r_{-t}(c)\| \, dt \, \,
\right)_{i \in I}
$$
converges to 0. \ \ \ \ \ (*)

\medskip

Choose $\vep > 0$.
By assumption, there exists a strictly positive number $M$ such that
$\|a_i\| \leq M$ for every $i \in I$ and $\|a\| \leq M$.

\begin{enumerate}
\item First we are going to approximate the integral over $\R$ by an integral over an
interval.

It is clear that there exists a natural number $m$ such that
$$  \int_{\R \setminus [-m,m]}  \exp(-r^2 t^2) \,  dt \leq \frac{\vep}{6 M
(  \|c\|+1)} \ .$$ Because $\|(a_i-a) \al^r_{-t}(c) \| \leq 2 M
\|c\|$ for every $i \in I$ and $t \in \R$, this implies that
\begin{eqnarray*}
\int \exp(-r^2 t^2) \, \| (a_i-a) \al^r_{-t}(c) \| \, dt
& = & \int_{-m}^m \exp(-r^2 t^2) \, \| (a_i-a) \al^r_{-t}(c) \| \, dt \\ &  & + \int_{\R
\setminus [-m,m]} \exp(-r^2 t^2) \, \| (a_i-a) \al^r_{-t}(c) \|
\, dt \\
& \leq &
\int_{-m}^m \exp(-r^2 t^2) \, \| (a_i-a) \al^r_{-t}(c) \| \, dt \ \ + \ \ \frac{\vep}{3}
\end{eqnarray*}
for every $i \in I$.

\item In a second step, we are going to approximate the integral over the interval by a sum.

By uniform continuity, there exists $\sde > 0$ such that we have for every $s,t \in
[-m,m]$ with $|s-t| \leq \sde$ that
$$|\ \exp(-r^2 t^2) \, \|(a_i-a) \al^r_{-t}(c) \|
- \exp(-r^2 s^2) \, \|(a_i-a) \al^r_{-s}(c) \| \ |
\leq  \frac{\vep}{6 m} \ .$$

Take a partition $r_0,\ldots\!,r_p$ for $[-m,m]$ such that $|r_j - r_{j-1}| \leq \sde$ for
$j=1,\ldots\!,p$. Because of the previous inequality, we get that
\begin{eqnarray*}
 \int_{-m}^m \exp(-r^2 t^2) \, \|(a_i - a) \al^r_{-t}(c) \| \, dt
& \leq &  \frac{\vep}{3} + \sum_{j=1}^p \,(r_j - r_{j-1}) \, \exp(-r^2 \, r_j^2) \,
\|(a_i - a) \al^r_{-r_j}(c)\|  \\
& \leq &  \frac{\vep}{3} + \sum_{j=1}^p (r_j - r_{j-1}) \, \|(a_i - a) \al^r_{-r_j}(c)\|
\end{eqnarray*}
for every $i \in I$. (The difference between upper- and lower sum for the integral under
this partition is smaller than $\frac{\vep}{3}$).

\item In the third step, we will use that $(a_i)_{i \in I}$ converges strictly to $a$ :

There exists $i_0 \in I$ such that we have for every $i \in I$ with $i \geq i_0$ and $j
\in \{1,\ldots\!,p\}$ that  \newline $\|(a_i - a) \al^r_{-r_j}(c)\| \leq \frac{\vep}{6 m}$.
This implies for every $i \in I$ with $i \geq i_0$ that
$$\sum_{j=1}^p (r_j - r_{j-1}) \,  \|(a_i - a) \al^r_{-r_j}(c)\|  \leq
\frac{\vep}{3} \ .$$
\end{enumerate}

Combining these three results, we see that
$$\int \exp(-r^2 t^2) \, \|(a_i - a) \al^r_{-t}(c)\| \, dt  \leq \vep$$
for every $i \in I$ with $i \geq i_0$.

\medskip

By remark (*), we see that $(a_i(r,z)\,b)_{i \in I}$ converges to $a(r,z)\,b$.

\medskip

Using $\al^l$, one proves in a similar way that $(b\,a_i(r,z))_{i \in I}$ converges to $b
\, a(r,z) $.
\end{demo}

The following proof is due to A. Van Daele and J. Verding.

\begin{proposition} \label{pa2.prop2}
Consider $z \in \C$. Then there exists a net $(u_k)_{k \in K}$ in $A$ consisting of
analytic elements for $\al$ and such that
\begin{itemize}
\item We have for every $k \in K$ that $\|u_k\| \leq 1$ and $\|\al_z(u_k)\| \leq 1$.
\item The nets $(u_k)_{k \in K}$ and $(\al_z(u_k))_{k \in K}$ converge strictly to 1.
\end{itemize}
\end{proposition}
\begin{demo}
Take an approximate identity $(e_q)_{q \in Q}$ for $A$. For every $q \in Q$ and $r > 0$,
we have the element $e_q(r)$ which is analytic with respect to $\al$ and satisfies
$\|e_q(r)\| \leq 1$ and $\|\al_z(e_q(r))\| \leq \exp(r^2 (\im z)^2)$.

By the previous result, we have for every $r > 0$ that $(e_q(r))_{q \in Q}$ and
$\bigl(\al_z(e_q(r))\bigr)_{q \in Q}$ are bounded nets which converge strictly to 1.

\medskip

Let us now define the set $$K = \{\, (F,n) \mid F \text{\ is a finite subset of  \ } A
\text{\ and \ } n \in \N  \, \} .$$ On $K$ we put an order such that
$$(F_1,n_1) \leq (F_2,n_2) \ \
\Leftrightarrow \ \ F_1 \subseteq F_2 \text{\ and \ } n_1 \leq n_2 $$
for every $(F_1,n_1),(F_2,n_2) \in K$.
In this way, $K$ becomes a directed set.

\medskip

Let us fix $k=(F,n) \in K$.

Firstly, there exist an element $r_k > 0$ such that $|\exp(-r_k^2 (\im z)^2) - 1 |
\leq \frac{1}{2 n} \frac{1}{\|x\|+1}$ for every $x \in F$.

Secondly, there exist $q_k \in Q$ such that
\begin{itemize}
\item We have that $\| e_{q_k}(r_k) \, x - x \| \leq \frac{1}{2 n}$
      for every $x \in F$.
\item We have that $\| x \, e_{q_k}(r_k) - x \| \leq \frac{1}{2 n}$
      for every $x \in F$.
\item We have that $\| \al_z(e_{q_k}(r_k)) \, x - x \|
      \leq \frac{1}{2 n}$ for every $x \in F$.
\item We have that $\| x \, \al_z(e_{q_k}(r_k)) - x \|
      \leq \frac{1}{2 n}$ for every $x \in F$.
\end{itemize}
Now we define $u_k = \exp(-r_k^2 (\im z)^2) \, e_{q_k}(r_k) \in A$.

It follows that $u_k$ is analytic with respect to $\al$ and  $\|u_k\| \leq 1$. We have
also that
$$
\|\al_z(u_k)\|  =  \exp(-r_k^2 (\im z)^2) \, \| \al_z(e_{q_k}(r_k))\|
\leq  \exp(r_k^2 (\im z)^2) \, \exp(-r_k^2 (\im z)) = 1
$$

\medskip

Now we prove that $(u_k)_{k \in K}$ converges strictly to 1.

Choose $x \in A$ and $\vep > 0$. Then there exists $n_0$ in $\N$ such that $\frac{1}{n_0}
\leq \vep$. Put $k_0 = (\{x\},n_0) \in K$.

Take $k=(F,n) \in K$ such that $k \geq k_0$, so $x$ belongs to $F$
and $n_0 \leq n$. Therefore,
\begin{eqnarray*}
& & \| u_k \, x - x \| =
\| \exp(-r_k^2 (\im z)^2) \,e_{q_k}(r_k) \, x - x \| \\
& & \spat \leq  \| \exp(-r_k^2 (\im z)^2) \, e_{q_k}(r_k) \, x
- \exp(-r_k^2 (\im z)^2) \, x \|
+ \| \exp(-r_k^2 (\im z)^2) \, x - x \| \\ & & \spat =  \exp(-r_k^2 (\im z)^2)
\, \| e_{q_k}(r_k) \, x - x \| + | \exp(-r_k^2 (\im z)^2) - 1 | \, \|x \| \\ & &
\spat \leq \| e_{q_k}(r_k) \, x - x \| + \frac{1}{2 n (\|x\|+1)} \, \|x\|
\leq  \frac{1}{2 n} + \frac{1}{2 n} = \frac{1}{n}  \leq  \frac{1}{n_0} \leq \vep \  .
\end{eqnarray*}
So we see that $(u_k \, x)_{k \in K}$ converges to $x$.
In a similar way, one proves that $(x \, u_k)_{k \in K}$ converges to $x$.

Consequently, the net $(u_k)_{k \in K}$ converges stricly to 1.

\medskip

Completely analogously, one proves that $\bigl(\al_z(u_k)\bigr)_{k \in K}$ converges
strictly to 1.
\end{demo}

\medskip

The following proposition follows immediately from propositions \ref{pa2.prop3} and
\ref{pa2.prop2} (applied to $\al^r$ instead of $\al$).

\begin{proposition} \label{pa2.prop5}
Consider $z \in \C$ and $a \in D(\oal_z)$.
Then there exists a net $(a_k)_{k \in K}$ in $D(\al_z)$ such that
\begin{enumerate}
\item We have that $\|a_k\| \leq \|a\|$ and $\|\al_z(a_k)\| \leq \|\al_z(a)\|$ for
every $k \in K$.
\item $(a_k)_{k \in K}$ converges strictly to $a$ and
$(\al_z(a_k))_{k \in K}$ converges strictly to $\al_z(a)$.
\end{enumerate}
\end{proposition}

\medskip

\begin{corollary} \label{pa2.cor1}
Consider $z \in \C$ and a core $C$ for $\al_z$. Let $a$ be an element in $D(\oal_z)$. Then
there exists a net $(a_k)_{k \in K}$ in $C$ such that
\begin{enumerate}
\item We have that $\|a_k\| \leq \|a\|$ and $\|\al_z(a_k)\| \leq \|\al_z(a)\|$ for
every $k \in K$.
\item $(a_k)_{k \in K}$ converges strictly to $a$ and
$(\al_z(a_k))_{k \in K}$ converges strictly to $\al_z(a)$.
\end{enumerate}
\end{corollary}
\begin{demo}
By the previous proposition, we get the existence of a net $(e_j)_{j \in J}$ in $D(\al_z)$
such that
\begin{enumerate}
\item We have that $\|e_j\| \leq \|a\|$ and $\|\al_z(e_j)\| \leq \|\al_z(a)\|$ for
every $j \in J$.
\item $(e_j)_{j \in J}$ converges strictly to $a$ and
$(\al_z(e_j))_{j \in J}$ converges strictly to $\al_z(a)$.
\end{enumerate}
Define $K = \N \times J$ and put on $K$ the product order. By lemma \ref{pa1.lem1}, there
exists for every $k=(n,j) \in K$ an element $a_k \in C$ such that
\begin{enumerate}
\item We have that $\|a_k\| \leq \|e_j\|$ and $\|\al_z(a_k)\| \leq \|\al_z(e_j)\|$
\item We have that $\|a_k - e_j\| \leq \frac{1}{n}$ and $\|\al_z(a_k) - \al_z(e_j)\|
\leq \frac{1}{n}$
\end{enumerate}
It is then easy to check that $(a_k)_{k \in K}$ satisfies the requirements of the
corollary.
\end{demo}

\bigskip\bigskip

In the last part of this section, we want to prove that $\oal_z$ is strictly closed. The
proof of this fact will be more involved than the proof of theorem \ref{pa1.th1}. We will
split it up in  several steps.

\begin{lemma} Consider $z \in \C$. Then the mapping $\al_z$ is strictly closable.
\end{lemma}
\begin{demo}
Choose a net $(a_i)_{i \in I}$ in $D(\al_z)$ and $b \in M(A)$ such that $(a_i)_{i \in I}$
converges strictly to $0$ and $(\al_z(a_i))_{i \in I}$ converges strictly to $b$.

Choose $c \in D(\al^r_z)$. We have for every $i \in I$ that $a_i \, c$ belongs to
$D(\al_z)$ and $\al_z(a_i \, c) = \al_z(a_i) \al^r_z(c)$. This implies immediately that
$(\al_z(a_i \, c))_{i \in I}$ converges to $b \, \al^r_z(c)$. It is also clear that $(a_i
\, c)_{i \in I}$ converges to 0. So the closedness of $\al_z$ implies that $b\,\al^r_z(c)
= 0$.

Because $\al^r_z$ has dense range, this implies that $b=0$.
\end{demo}

\bigskip

We will now introduce a temporary notation.

\begin{notation}
Consider $z \in \C$. Then $\tau_z$ will denote the strict closure of $\al_z$ in $M(A)$.
\end{notation}

We will ultimately prove that $\tau_z = \oal_z$.

\medskip

\begin{result} \label{pa2.res2}
Let $z$ be an element in $\,\C$ and $a \in D(\tau_z)$. Then we have the following
properties :
\begin{itemize}
\item We have for every $b \in D(\al^l_z)$ that $b \, a \in D(\al_z)$ and $\al_z(b \, a) =
\al^l_z(b) \tau_z(a)$.
\item We have for every $b \in D(\al^r_z)$ that $a \, b \in D(\al_z)$ and $\al_z(a \, b) =
\tau_z(a) \al^r_z(b)$.
\end{itemize}
\end{result}
\begin{demo} Choose $b \in D(\al^r_z)$.
There exist a net $(a_i)_{i \in I}$ in $D(\al_z)$ such that $(a_i)_{i \in I}$ converges
strictly to $a$ and $(\al_z(a_i))_{i \in I}$ converges strictly to $\tau_z(a)$.

We have for every $i \in I$ that $a_i \, b$ belongs to $D(\al_z)$ and $\al_z(a_i \, b) =
\al_z(a_i) \al^r_z(b)$. This implies that $(\al_z(a_i \, b))_{i \in I}$ converges to
$\tau_z(a) \al^r_z(b)$. Because we have also that $(a_i \, b)_{i \in I}$ converges to $a
\, b$, the closedness of $\al_z$ implies that $a \, b$ belongs to $D(\al_z)$ and
$\al_z(a \, b) = \tau_z(a) \al^r_z(b)$.

\medskip

The result concerning $\al^l$ is proven in a similar way.
\end{demo}

There is a converse of this result :

\begin{result} \label{pa2.res3}
Consider $z \in \C$ and $a \in M(A)$ such that there exists an element $c \in M(A)$ such
that we have for every $b$ in $D(\al^r_z)$ that $a \, b$ belongs to $D(\al_z)$ and
$\al_z(a \, b) = c \, \al^r_z(b)$. Then $a$ belongs to $D(\tau_z)$ and $\tau_z(a)=c$.
\end{result}
\begin{demo}
We know by proposition \ref{pa2.prop2} (applied to $\al^r$ instead of $\al$) that there
exists a net $(u_k)_{k \in K}$ in $D(\al^r_z)$ such that $(u_k)_{k \in K}$ converges
strictly to 1 and $(\al^r_z(u_k))_{k \in K}$ converges strictly to 1.

We have by assumption for every $k \in K$ that $a \, u_k$ belongs to $D(\al_z)$ and
$\al_z(a \, u_k) = c \, \al^r_z(u_k)$. Consequently, we have that $(\al_z(a \, u_k))_{k
\in K}$ converges strictly to $c$. It is also clear that $(a \, u_k)_{k \in K}$ converges
strictly to $a$. By definition, we find that $a$ belongs to $D(\tau_z)$ and $\tau_z(a) =
c$.
\end{demo}

There is of course an analogous result concerning $\al^l$.

\medskip

Using proposition \ref{pa2.prop2}  (applied to $\al^r$ instead of $\al$) and result
\ref{pa2.res2}, we get also easily the following lemma.

\begin{lemma} \label{pa2.lem2}
Let $z$ be a complex number and $a \in D(\tau_z)$. Then there exists a net $(a_k)_{k \in
K}$ in $D(\al_z)$ such that
\begin{enumerate}
\item We have that $\|a_k\| \leq \|a\|$ and $\|\al_z(a_k)\| \leq \|\tau_z(a)\|$ for every
$k \in K$.
\item $(a_k)_{k \in K}$ converges strictly to $a$ and $(\al_z(a_k))_{k \in K}$ converges
strictly to $\tau_z(a)$.
\end{enumerate}
\end{lemma}

We also have the following result.

\begin{lemma} \label{pa2.lem3}
Let $y,z$ be complex numbers such that $y$ belongs to $S(z)$. Consider $a \in D(\tau_z)$
and a net $(a_k)_{k \in K}$ in $D(\al_z)$ such that $(a_k)_{k \in K}$ converges strictly
to $a$ and $(\al_z(a_k))_{k \in K}$ converges strictly to $\tau_z(a)$. Suppose moreover
that there exists a positive number $M$ such that $\|a_k\| , \|\al_z(a_k)\| \leq M$ for
every $k \in K$. Then $a$ belongs to $D(\tau_y)$, we have the inequality   $\|\tau_y(a)\|
\leq M$ and $(\al_y(a_k))_{k \in K}$ converges strictly to $\tau_y(a)$.
\end{lemma}
\begin{demo}
By corollary \ref{pa1.cor1}, we know for every $k \in K$ that $a_k$ belongs to $D(\al_y)$
and $\|\al_y(a_k)\| \leq \max\{\|a_k\|, \|\al_z(a_k)\|\}$ $\leq M$. Choose $b \in
D(\al^r_z)$. Then we have for ever $k,l \in K$ that
\begin{eqnarray*}
& & \|\al_y(a_k) \al^r_y(b) - \al_y(a_l) \al^r_y(b) \|  = \|\al_y(a_k b - a_l b)\| \\ & &
\spat \leq \max\{\|a_k b - a_l b\| , \|\al_z(a_k b - a_l b)\| \}
= \max\{\|a_k b - a_l b\| , \|\al_z(a_k) \al^r_z(b) - \al_z(a_l) \al^r_z(b)\|\}
\end{eqnarray*}
This implies that $\bigl(\al_y(a_k) \al^r_y(b)\bigr)_{k \in K}$ is Cauchy and hence
convergent.

We have already seen that $(\al_y(a_k))_{k \in K}$ is bounded. Combining this with the
previous fact gives us that $(\al_y(a_k) c)_{k \in K}$ is convergent for every $c \in A$.

Using $\al^l$ in stead of $\al^r$, one proves similarly that $(c \, \al_y(a_k))_{k \in K}$
is convergent for every $c \in A$.

These two facts imply that $(\al_y(a_k))_{k \in K}$ is strictly convergent in $M(A)$. By
definition, this implies that $a$ belongs to $D(\tau_y)$ and that $(\al_y(a_k))_{k \in K}$
converges strictly to $\tau_y(a)$. Because $\|\al_y(a_k)\| \leq M$ for every  $k \in K$, we
have also immediately that $\|\tau_y(a)\| \leq M$.
\end{demo}

Combining these two lemmas, we get easily the following one.

\begin{lemma} \label{pa2.lem1}
Let $y,z$ be complex numbers such that $y \in S(z)$. Then $D(\tau_z) \subseteq D(\tau_y)$
and we have for every $a \in D(\tau_z)$ that $\|\tau_y(a)\| \leq \max\{\|a\|,\|\tau_z(a)\|\}$
\end{lemma}

\medskip

\begin{proposition} \label{pa2.prop4}
Consider $z \in \C$ and $a \in D(\tau_z)$. Then the function
$S(z) \rightarrow M(A) : y \mapsto \tau_y(a)$ is strictly continuous.
\end{proposition}
\begin{demo} Take $b \in A$.

Choose $y \in S(z)$. Take $\vep > 0$. By the previous lemma, there exists a strictly
positive number $M$ such that $\|\tau_u(a)\| \leq M$ for every $u \in S(z)$. Now we have :
\begin{itemize}
\item There exists $c \in D(\al^r_z)$ such that $\|\al^r_y(c) - b\| \leq \frac{\vep}{4M}$.
\item Because $c \in D(\al^r_z)$,
there exists $\sde_1 > 0$ such that $\| \al^r_u(c) - \al^r_y(c) \| \leq
\frac{\vep}{4M}$ for every $u \in S(z)$ with $|u-y| \leq \sde_1$.
\item We know that $a \, c$ belongs to $D(\al_z)$ by result \ref{pa2.res2}. This
implies the existence of $\sde_2 > 0$ such that $\| \al_u(a \, c) - \al_y(a \, c) \| \leq
\frac{\vep}{4}$ for every $u \in S(z)$ with $|u-y| \leq \sde_2$.
\end{itemize}
Put $\sde = \min\{\sde_1,\sde_2\}$. Then we have for every $u \in S(z)$ with $|u-y|
\leq \sde$ that
\begin{eqnarray*}
\| \tau_u(a) b - \tau_y(a) b \| & \leq & \| \tau_u(a) b - \tau_u(a) \al^r_y(c)\| +
\| \tau_u(a) \al^r_y(c) - \tau_u(a) \al^r_u(c) \| \\
& & + \| \tau_u(a) \al^r_u(c) - \tau_y(a) \al^r_y(c) \| + \| \tau_y(a)  \al^r_y(c) -
\tau_y(a) b \| \\
& \leq & \|\tau_u(a)\| \, \| b - \al^r_y(c) \| + \|\tau_u(a)\| \, \|\al^r_y(c) -
\al^r_u(c)\| \\ & & + \| \al_u(a\,c) - \al_y(a\,c)\| + \| \tau_y(a)\| \, \|\al^r_y(c)
-  b \| \\ &
\leq & M \frac{\vep}{4 M} + M \frac{\vep}{4M} + \frac{\vep}{4} + M \frac{\vep}{4M} = \vep
\ .
\end{eqnarray*}
So we have proven that the mapping $S(z) \rightarrow A : z \mapsto
\tau_z(a) b$ is continuous.

\medskip

In a similar way, one proves that the mapping $S(z) \rightarrow A : z \mapsto b \,
\tau_z(a)$ is continuous.
\end{demo}

\bigskip

In order to  prove analyticity on the interior, we will need first a small
estimation.

\begin{lemma}
Consider $z \in \C$, $y \in S(z)^0$ and put $\sde = \min\{|\im y|,|\im y-
\im z\,|\} > 0$. Let $a$ be an element in $D(\al_z)$.
Then we have for every $u \in \C$ with $|u-y| \leq \frac{\sde}{4}$ that
$$\| \al_u(a) - \al_y(a) \| \leq \frac{8}{\sde} \max\{\|a\|,\|\al_z(a)\|\} \,\, |u-y| \ .$$
\end{lemma}
\begin{demo}
Define the function $f$ from $S(z)^0$ into $A$ such that $f(u) = \al_u(a)$ for every $u
\in S(z)^0$. Then $f$ is analytic on $S(z)^0$. Define the curve $\gam : [0,2 \pi]
\rightarrow S(z) : \th \mapsto \frac{\sde}{2} \, e^{i   \th}$ in $S(z)^0$.

Choose $v \in \C$ with $|v-y| \leq \frac{\sde}{4}$.
By the Cauchy theorem, we have that
$$f'(v)=  \frac{1}{2 \pi i} \int_\gam \frac{f(p)}{(p-v)^2} \, dp \ . $$
We also know, using corollary \ref{pa1.cor1} , that
$$ \frac{\|f(p)\|}{|p-v|^2} = \frac{\|\al_p(a)\|}{|p-v|^2} \leq
\max\{\|a\|,\|\al_z(a)\|\} \, \frac{16}{\sde^2} $$
for every $p \in \bld \gam$. Combining these two facts, we see that $$\|f'(v)\|
\leq \frac{1}{2 \pi} \, \max\{\|a\|,\|\al_z(a)\|\} \, \frac{16}{\sde^2} \,  \pi \sde =
\frac{8}{\sde} \, \max\{\|a\|,\|\al_z(a)\|\} \ .$$

\medskip

Using a mean value theorem, we get for every $u \in \C$ with $|u-y| \leq \frac{\sde}{4}$
that
$$\|\al_u(a) - \al_y(a)\| = \|f(u) - f(y)\| \leq
\frac{8}{\sde} \, \max\{\|a\|,\|\al_z(a)\|\} \,\, |u-y| .$$
\end{demo}

\medskip

We want to generalize this lemma to $\tau$ :

\begin{lemma}
Consider $z \in \C$, $y \in S(z)^0$ and put $\sde =  \min\{|\im y|,|\im y-
\im z\,|\} > 0$. Let $a$ be an element in $D(\tau_z)$. Then we have for every $u \in
\C$ with $|u-y| \leq \frac{\sde}{4}$ that
$$\| \tau_u(a) - \tau_y(a) \| \leq \frac{8}{\sde} \max\{\|a\|,\|\tau_z(a)\|\}
\,\, |u-y| \ .$$
\end{lemma}
\begin{demo}
By lemma \ref{pa2.lem2}, there exists a net $(a_k)_{k \in K}$ in $D(\al_z)$ such that
\begin{enumerate}
\item We have that $\|a_k\| \leq \|a\|$ and $\|\al_z(a_k)\| \leq \|\tau_z(a)\|$ for every
$k \in K$.
\item $(a_k)_{k \in K}$ converges strictly to $a$ and
$(\al_z(a_k))_{k \in K}$ converges strictly to $\tau_z(a)$.
\end{enumerate}

Choose $u \in \C$ with $|u - y| \leq \frac{\sde}{4}$.
By lemma \ref{pa2.lem3}, we know that $(\al_u(a_k) - \al_y(a_k))_{k \in K}$ converges
strictly to $\tau_u(a) - \tau_y(a)$.
By the previous lemma, we have for every $k \in K$
that
\begin{eqnarray*}
\|\al_u(a_k) - \al_y(a_k)\| & \leq & \frac{8}{\sde}  \max\{\|a_k\|,\|\al_z(a_k)\|\}
\,\, |u-y| \\
& \leq & \frac{8}{\sde} \max\{\|a\|,\|\tau_z(a)\|\} \,\, |u-y| \ .
\end{eqnarray*}

Therefore, we get that $\|\tau_u(a) - \tau_y(a)\| \leq
\frac{8}{\sde} \max\{\|a\|,\|\tau_z(a)\|\} \,\, |u-y| $.
\end{demo}

\begin{proposition}
Consider $z \in \C$ and $a \in D(\tau_z)$. Then the function
$S(z)^0 \rightarrow M(A) : y \mapsto \tau_y(a)$ is analytic.
\end{proposition}
\begin{demo}
Choose $b \in A$.

Take $y \in S(z)^0$. Choose $\vep > 0$. By lemma \ref{pa2.lem1}, we have the existence of
a strictly positive number $N$ such that $\| \tau_u(a) \| \leq N$ for every $u \in S(z)$.

Put $\sde_1 = \frac{1}{4} \, \min\{|\im y|,|\im y- \im z\,|\} > 0$. Using the previous
lemma, we get the existence of  a strictly positive number $M$ such that $\|
\tau_u(a) - \tau_y(a) \| \leq M \, |u-y|$ for every $u \in S(z)^0$ with $|u-y| \leq
\sde_1$ \ \ \ \ (*).

Furthermore, there exists $c \in D(\al^r_z)$ such that $\|b - \al^r_y(c) \| \leq
\frac{\vep}{4 M}$. We know that $a \, c$ belongs to $D(\al_z)$ by result \ref{pa2.res2}.

Define the functions $f,g$ from $S(z)^0$ into $A$ such that $f(u)= \al_u(a \, c)$ and
$g(u) = \al^r_u(c)$ for every $u \in S(z)^0$. Then $f$ and $g$ are analytic on $S(z)^0$.
We have moreover the following properties :
\begin{itemize}
\item There exists $\sde_2 > 0$ such that
$$\| \frac{\al^r_u(c) - \al^r_y(c)}{u-y} - g'(y) \| \leq \frac{\vep}{12 N} $$
for every $u \in S(z)^0$ with $0 < |u-y| \leq \sde_2$.
\item There exists $\sde_3 > 0$ such that
$$\| \frac{\al_u(a\,c) - \al_y(a\,c)}{u-y} - f'(y) \| \leq \frac{\vep}{12}  $$
for every $u \in S(z)^0$ with $0 < |u-y| \leq \sde_3$.
\item By proposition \ref{pa2.prop4}, there exists $\sde_4 > 0$ such that
$$\| \tau_u(a) g'(y) - \tau_y(a) g'(y) \| \leq \frac{\vep}{12}$$
for every $u \in S(z)^0$ with $|u-y| \leq \sde_4$.
\end{itemize}
Put $\sde = \min\{\sde_1,\sde_2,\sde_3,\sde_4\}$. The 3 above inequalities imply for all
$u \in S(z)^0$ with $0 < |u-y| \leq \sde$ that
\begin{eqnarray*}
& & \| \frac{\tau_u(a) \al^r_y(c) - \tau_y(a) \al^r_y(c)}{u-y} - f'(y) + \tau_y(a) g'(y) \|
\\ & & \spat =  \| \frac{\tau_u(a) \al^r_y(c) - \tau_u(a) \al^r_u(c) + \tau_u(a)
\al^r_u(c) -
\tau_y(a) \al^r_y(c)}{u-y} \\
& & \spat \ \ \ + \tau_u(a) g'(y) - \tau_u(a) g'(y) - f'(y) + \tau_y(a) g'(y) \|
\\ & & \spat \leq \| \tau_u(a) \, \left( \frac{\al^r_y(c) - \al^r_u(c)}{u - y} + g'(y)
\right) \| + \| \frac{\tau_u(a) \al^r_u(c) - \tau_y(a) \al^r_y(c)}{u-y} - f'(y)\|
\\ & & \spat\ \ \ + \| - \tau_u(a) g'(y) + \tau_y(a) g'(y) \| \\
& & \spat \leq \|\tau_u(a)\|
\, \|\frac{\al^r_y(c) - \al^r_u(c)}{u - y} + g'(y)\| + \| \frac{\al_u(a\,c) -
\al_y(a\,c)}{u-y}
- f'(y) \| + \frac{\vep}{12} \\
& & \spat \leq N \frac{\vep}{12 N} + \frac{\vep}{12} + \frac{\vep}{12} = \frac{\vep}{4} \ .
\end{eqnarray*}
This implies for every $u_1,u_2 \in S(z)^0$ with $0 < |u_1 - y| \leq \sde$ and   $0 < |u_2
- y| \leq \sde$ that
\begin{eqnarray*}
& & \| \frac{\tau_{u_1}(a) b - \tau_y(a) b}{u_1 - y} - \frac{\tau_{u_2}(a) b - \tau_y(a) b}
{u_2 - y} \| \\
& & \spat \leq  \| \frac{\tau_{u_1}(a) b - \tau_y(a) b}{u_1 - y}
- \frac{\tau_{u_1}(a) \al^r_y(c)  - \tau_y(a) \al^r_y(c)}{u_1 - y} \| \\
& & \spat \ \ \  + \| \frac{\tau_{u_1}(a) \al^r_y(c) - \tau_y(a) \al^r_y(c)}{u_1 - y}
- f'(y) + \tau_y(a) g'(y) \| \\
& & \spat\ \ \ + \| f'(y) - \tau_y(a) g'(y) - \frac{\tau_{u_2}(a) \al^r_y(c) -
\tau_y(a) \al^r_y(c)}{u_2 - y} \| \\
& & \spat\ \ \ + \| \frac{\tau_{u_2}(a) \al^r_y(c) - \tau_y(a) \al^r_y(c)}{u_2 - y} -
\frac{\tau_{u_2}(a) b - \tau_y(a) b}{u_2 - y} \| \\
& & \spat \leq \| \frac{\tau_{u_1}(a) - \tau_y(a)}{u_1 - y} \| \, \|b - \al^r_y(c)\| +
\frac{\vep}{4} + \frac{\vep}{4}
+ \| \frac{\tau_{u_2}(a) - \tau_y(a)}{u_2 - y} \| \, \|\al^r_y(c) - b\|  \\
& & \spat \leq M \frac{\vep}{4 M} + \frac{\vep}{4} + \frac{\vep}{4} + M \frac{\vep}{4 M}
= \vep \ ,
\end{eqnarray*}
where (*) was used in the last inequality.

\medskip

From this all, we infer that the function $S(z)^0 \rightarrow A : y \mapsto \tau_y(a) b$ is
analytic. This implies that the function $S(z)^0 \rightarrow M(A) : y \mapsto \tau_y(a)$ is
analytic.
\end{demo}

\medskip

So we arrive at the following conclusion.

\begin{proposition}
We have for every $z \in \C$ that $\tau_z=\oal_z$.
\end{proposition}

The previous results imply easily that $\tau_z \subseteq \oal_z$.
The other one follows from proposition \ref{pa2.prop5}.

\bigskip

Therefore, we have proven the following theorem.

\begin{theorem} \label{pa2.th1}
Consider $z \in \C$. Then the mapping $\al_z$ is strictly closed and
has $D(\al_z)$ as strict core.
\end{theorem}

\bigskip

The next proposition depends essentially on the strict closedness of $\oal_z$
(see the proof of lemma \ref{pa2.res3}).

\begin{proposition}
Consider elements $a,b \in M(A)$ and $z \in \C$. Then
\begin{enumerate}
\item If we have for every $c \in D(\al^l_z)$ that $c\,a$ belongs to $D(\al_z)$ and
$\al_z(c \, a) = \al^l_z(c) \, b$, then $a$ belongs to $D(\oal_z)$
and $\al_z(a) = b$.
\item If we have for every $c \in D(\al^r_z)$ that $a \, c$ belongs to $D(\al_z)$ and
$\al_z(a \, c) = b \, \al^r_z(c)$, then $a$ belongs to $D(\oal_z)$
and $\al_z(a) = b$.
\end{enumerate}
\end{proposition}

\bigskip

Because $\al^l$ and $\al^r$ satisfy the same conditions as $\al$, we can define for every
$z \in \C$ the strictly closed operators $\overline{\al^l}_z$ and $\overline{\al^r}_z$
from within $M(A)$ into $M(A)$. Again we put $\al^l_z(a) = \overline{\al^l}_z(a)$ for
every $a \in D(\overline{\al^l}_z)$ and $\al^r_z(a) = \overline{\al^r}_z(a)$ for every $a
\in D(\overline{\al^r}_z)$.

\medskip

Then it is easy to prove the following result :

\begin{proposition}
Consider a complex number $z$ and an element $a \in D(\overline{\al}_z)$. Then we have the
following properties :
\begin{enumerate}
\item We have for every $b \in D(\overline{\al^l}_z)$  that $b\,a$ belongs to $D(\oal_z)$
and $\al_z(b\,a) = \al^l_z(b) \al_z(a)$.
\item We have for every $b \in D(\overline{\al^r}_z)$  that $a\,b$ belongs to $D(\oal_z)$
and $\al_z(a\,b) = \al_z(b) \al^r_z(b)$.
\end{enumerate}
\end{proposition}

\bigskip

The proof of the next proposition is completely similar to te proof of proposition
\ref{pa1.prop3}.

\begin{proposition} \label{pa2.prop8}
Consider a subset $K$ of $A^*$ such that $K$ is separating for $M(A)$ and $\th  \, \al_t
\in K$ for every $\th \in K$ and $t \in \R$. Let $z \in \C$  and $a,b \in M(A)$  such that
there exists for every $\th \in K$ a function $f_\th$ from $S(z)$ in $\C$ such that
\begin{enumerate}
\item $f_\th$ is bounded and continuous on $S(z)$
\item $f_\th$ is analytic on $S(z)^0$
\item We have for every $t \in \R$ that $f_\th(t) = \th(\al_t(a))$
\item We have that $f_\th(z) = \th(b)$.
\end{enumerate}
Then $a$ belongs to $D(\oal_z)$ and $\al_z(a) = b$.
\end{proposition}

\bigskip

\section{Commuting one-parameter representations}

Consider a Banach space $E$ and two norm continuous one-parameter representations $\al$
and $\be$ on $E$ which commute : $\al_s \, \be_t = \be_t \, \al_s$ for every $s,t \in \R$.

We define the norm continuous one-parameter representation $\gamma$ on $E$ such that
$\gamma_t = \al_t \, \be_t$ for every $t \in \R$.

It is clear that $\al$, $\be$ and $\gam$ mutually commute.

\medskip

Using merely the definition of the analytic continuation of a one-parameter
representations, it is not so difficult to prove for every $z \in \C$ and $t \in \R$ that
$\al_t \, \be_z = \be_z \, \al_t$ and that $\be_t \, \al_z = \al_z \, \be_t$

The same relations hold of course between $\al$,$\gam$ and $\be$,$\gam$.

\medskip

Consider $z \in \C$. In the following part, we will prove that $\al_z \, \be_z$ and
$\be_z \, \al_z$  are closable and that their closure is equal to $\gamma_z$.

\bigskip

It is easy to prove that the function $\R \times \R \rightarrow E : (s,t) \mapsto
\al_s(\be_t(a))$ is continuous.

\medskip

In the first part of this section we will use the following notation. We define for every
$a \in E$ and $n \in \N$ the  element
$$a(n) = \frac{n^2}{\pi} \, \int \int \exp(-n^2(s^2 + t^2)) \, \al_s(\be_t(a)) \, ds \,
dt \ .$$

\medskip

\begin{lemma} \label{pa4.lem1}
Consider $a \in E$, $z \in \C$ and $n \in \N$. Then we have the following properties :
\begin{enumerate}
\item The element $a(n)$ belongs to $D(\al_z \, \be_z) \cap D(\be_z \, \al_z)$
and $$(\al_z \, \be_z)(a(n)) = (\be_z \, \al_z)(a(n)) =  \frac{n^2}{\pi} \int
\exp\bigl(-n^2 ( (s-z)^2 + (t-z)^2)\bigr) \, \al_s(\be_t(a)) \, ds \, dt \ .$$
\item If $a$ belongs to $D(\al_z \, \be_z)$, then
$ (\al_z \, \be_z)(a(n)) = \bigl((\al_z \, \be_z)(a)\bigr)(n)  $.
\item If $a$ belongs to $D(\be_z \, \al_z)$, then
$(\be_z \, \al_z)(a(n)) =  \bigl((\be_z \, \al_z)(a)\bigr)(n) $.
\end{enumerate}
\end{lemma}
\begin{demo}
\begin{enumerate}
\item
We have that
$$a(n) = \frac{n^2}{\pi} \int \exp(-n^2 t^2) \, \be_t \left(
\int \exp(-n^2 s^2) \, \al_s(a)\, ds \right) \, dt \ .$$
This implies that $a(n)$ belongs to $D(\be_z)$ and
\begin{eqnarray*}
 \be_z (a(n))
& = & \frac{n^2}{\pi} \int \exp(-n^2 (t-z)^2) \, \be_t \left(
\int \exp(-n^2 s^2) \, \al_s(a)\, ds \right) \, dt  \\
& = & \frac{n^2}{\pi} \int \exp(-n^2 s^2) \, \al_s \left(
\int \exp(-n^2 (t-z)^2) \, \be_t(a)\, dt \right) \, ds
\end{eqnarray*}
and this implies that $\be_z(a(n))$ belongs to $D(\al_z)$  and that
\begin{eqnarray*}
\al_z(\be_z(a(n)))
& = & \frac{n^2}{\pi} \int  \exp(-n^2 (s-z)^2) \, \al_s \left(
\int \exp(-n^2 (t-z)^2) \, \be_t(a)\, dt \right) \, ds   \hspace{1cm} \text{(*)} \\
& = &  \frac{n^2}{\pi} \int \exp\bigl(-n^2 ( (s-z)^2 + (t-z)^2)\bigr) \, \al_s(\be_t(a))
\, ds \, dt \ .
\end{eqnarray*}

\item Because $a$ belongs to $D(\be_z)$, equation (*) implies that
$$\al_z(\be_z(a(n))) = \frac{n^2}{\pi} \int  \exp(-n^2 (s-z)^2) \, \al_s \left(
\int \exp(-n^2 t^2) \, \be_t(\be_z(a))\, dt \right) \, ds  \ . $$
So
$$\al_z(\be_z(a(n))) = \frac{n^2}{\pi} \int \exp(-n^2 t^2) \, \be_t \left(
\int \exp(- n^2 (s-z)^2) \, \al_s( \be_z(a) ) \, ds \right) \, dt .$$
Because $\be_z(a)$ belongs to $D(\al_z)$, this last equation implies that
\begin{eqnarray*}
\al_z(\be_z(a(n))) & = & \frac{n^2}{\pi} \int \exp(-n^2 t^2) \, \be_t \left(
\int \exp(- n^2 s^2) \, \al_s\bigl(\,\al_z(\be_z(a))\,\bigr) \, ds \right) \, dt \\
& = &  \frac{n^2}{\pi} \int \int \exp(-n^2(s^2 + t^2)) \,
\al_s\bigl(\be_t\bigl(\,(\al_z\,\be_z)(a)\,\bigr)\bigr) \, ds \, dt \ .
\end{eqnarray*}

\end{enumerate}

The results concerning $\be_z \, \al_z$ are proven in a similar way.
\end{demo}

\medskip

We can now use this lemma to prove the following result.

\begin{result}
Consider $z \in \C$. Then $\al_z \, \be_z$ and $\be_z \, \al_z$ are closable.
\end{result}
\begin{demo}
Choose a sequence $(a_j)_{j=1}^\infty$ in $D(\al_z \, \be_z)$ and $b \in E$ such that
$(a_j)_{j=1}^\infty$ converges to $0$ and $\bigl((\al_z \, \be_z)(a_j)\bigr)_{j=1}^\infty$
converges to $b$

\medskip

Take $m \in \N$.

By the  first statement of the previous lemma, we have for every $j \in \N$ that $a_j(m)$
belongs to $D(\al_z \, \be_z)$ and that the net $\bigl(\,(\al_z
\, \be_z)(a_j(m))\,\bigl)_{j=1}^\infty$ converges to 0.

By the second statement of the previous lemma, we get that the net $\bigl(\,(\al_z
\, \be_z)(a_j(m))\,\bigl)_{j=1}^\infty$ converges to $b(m)$. So we must have that $b(m) = 0$.

\vspace{1mm}

Because $(b(n))_{n=1}^\infty$ converges to $b$, we see that $b=0$.

\medskip

One proves in the same way that $\be_z \,\al_z$ is closable.
\end{demo}

\medskip

So we can introduce the following notation  :

\begin{notation}
Consider $z \in \C$. Then we define $\al_z \comp \be_z$ as the closure of $\al_z \,
\be_z$ and we define $\be_z \comp \al_z$ as the closure of $\be_z \,
\al_z$
\end{notation}

\medskip

The following lemma follows rather easily from lemma \ref{pa4.lem1}.

\begin{lemma} \label{pa4.lem2}
Consider $z \in \C$ and $n \in \N$. Then we have the following properties
\begin{itemize}
\item We have for every $a \in D(\al_z \comp \be_z)$ that
$(\al_z \comp \be_z)(a(n)) = \bigl((\al_z \comp \be_z)(a)\bigl)(n)$ .
\item We have for every $a \in D(\be_z \comp \al_z)$ that
$(\be_z \comp \al_z)(a(n)) = \bigl((\be_z \comp \al_z)(a)\bigl)(n)$ .
\end{itemize}
\end{lemma}

\medskip

Then we get the following proposition.

\begin{proposition} \label{pa4.prop2}
Consider $z \in \C$. Then $\al_z \comp \be_z = \be_z \comp \al_z$
\end{proposition}
\begin{demo}
Choose $a \in D(\al_z \comp \be_z)$. By lemma \ref{pa4.lem1}, we have for every $n
\in \N$ that $a(n)$ belongs to $D(\al_z \comp \be_z) \cap D(\be_z \comp \al_z)$ and
that $(\al_z \comp \be_z)(a(n)) = (\be_z \comp \al_z)(a(n))$.

\medskip

It is clear that $\bigl(a(n)\bigr)_{n=1}^\infty$ converges to $a$.

By the first statement of the previous lemma, we see that $\bigl((\al_z \comp
\be_z)(a(n))\bigl)_{n=1}^\infty$ converges to $(\al_z \comp \be_z)(a)$.

So $\bigl((\be_z \comp \al_z)(a(n))\bigl)_{n=1}^\infty$ converges also to $(\al_z \, .
\, \be_z)(a)$.

Hence the closedness of $\be_z \comp \al_z$ implies that $a$ belongs to $D(\be_z \comp
\al_z)$ and that $(\be_z \comp \al_z)(a) = (\al_z \comp
\be_z)(a)$.

\medskip

So we have proven that $\al_z \comp \be_z \subseteq \be_z \comp \al_z$. One proves in
a similar way that $\be_z \comp \al_z \subseteq \al_z \comp \be_z$.
\end{demo}

\medskip

The following result follows easily from lemmas \ref{pa4.lem1} and \ref{pa4.lem2}.

\begin{result}
Consider a complex number $z$. Then $D(\al_z \, \be_z) \cap D(\be_z \, \al_z)$ is a
core for $\al_z \comp \be_z$
\end{result}

\medskip

\begin{remark} \rm
The foregoing results (and their proofs) remain valid if we would look at compositions
$\al_y \, \be_z$ and $\be_z \, \al_y$. But we will be only interested in the case where $y
= z$.
\end{remark}

\bigskip

Now we are going to prove that $\gam_z = \al_z \comp \be_z = \be_z \comp \al_z$.

\begin{lemma}
Consider $a \in A$, $z \in \C$ and $n \in \N$. Then we have the following properties :
\begin{enumerate}
\item The element $a(n)$ belongs to $D(\al_z \, \be_z) \cap D(\be_z \, \al_z) \cap
D(\gamma_z) $ and $(\al_z \, \be_z)(a(n)) = (\be_z \, \al_z)(a(n)) =  \gamma_z(a(n))$.
\item If $a$ belongs to $D(\gamma_z)$, then
$ \gamma_z(a(n)) = (\gamma_z(a))(n)$ .
\end{enumerate}
\end{lemma}
\begin{demo}
\begin{enumerate}
\item Because the function $$
\C \rightarrow A : y \mapsto  \frac{n^2}{\pi} \int \int \exp\bigl(-n^2 (
(s-y)^2 + (t-y)^2) \bigr) \, \al_s(\be_t(a)) \, ds \, dt \ $$ is analytic
and
\begin{eqnarray*}
\gam_u(a(n))
& = &  \frac{n^2}{\pi} \int \int \exp\bigl(-n^2 ( s^2 + t^2) \bigr) \, \al_s(\be_t(\,
\gamma_u(a)\, )) \, ds \, dt  \\
& = & \frac{n^2}{\pi} \int \int \exp\bigl(-n^2 ( (s-u)^2 + (t-u)^2) \bigr)
\, \al_s(\be_t(a))
\, ds \, dt \\
\end{eqnarray*}
for every $u \in \R$, we have that $a(n)$ belongs to $D(\gamma_z)$ and
$$ \gamma_z(a(n))=  \frac{n^2}{\pi} \int \int \exp\bigl(-n^2 ( (s-z)^2 + (t-z)^2 )\bigr) \,
\al_s(\be_t(a)) \, ds \, dt \ .$$

So we see that $(\al_z \,\be_z)(a(n)) = (\be_z \,\al_z)(a(n)) = \gamma_z(a(n))$.

\item We know that $\al_s \gamma_z = \gamma_z \al_s$ for $s \in \R$
and that $\be_t \gamma_z = \gamma_z \be_t$ for every $t \in \R$.

So we have for every $s,t \in \R$ that $\al_s(\be_t(a))$ belongs to $D(\gamma_z)$ and
$\gamma_z\bigl(\al_s(\be_t(a))\bigr) = \al_s(\be_t(\, \gamma_z(a)\, ))$. This implies that
the function
$$ \R \times \R \rightarrow A : (s,t) \mapsto \exp(-n^2(s^2+t^2)) \,
\gamma_z\bigl( \al_s(\be_t(a)) \bigr)$$ is integrable.

Consequently, lemma \ref{o.lem1} and the closedness of
$\gamma_z$ imply that $a(n)$ belongs to $D(\gamma_z)$ and
\begin{eqnarray*}
\gamma_z(a(n))
& = & \frac{n^2}{\pi} \int \int \exp(-n^2(s^2+t^2)) \,
\gamma_z\bigl( \al_s(\be_t(a)) \bigr) \, ds \, dt \\
& = &
\frac{n^2}{\pi} \int \int \exp(-n^2(s^2+t^2)) \,
\al_s(\be_t(\, \gamma_z(a)\, )) \, ds \, dt
= (\gam_z(a))(n)\ .
\end{eqnarray*}
\end{enumerate}
\end{demo}

\medskip

Using this lemma, the proof of the next proposition is completely similar to the proof of
proposition \ref{pa4.prop2}.

\begin{proposition} \label{pa4.prop1}
Consider $z \in \C$. Then $\gamma_z = \al_z \comp \be_z  = \be_z \comp \al_z$
\end{proposition}

\bigskip\bigskip

For the rest of this section, we consider a \cst-algebra $A$ and strongly continuous
one-parameter representations $\al$,$\be$ on $A$ which commute and such that $\al$ and
$\be$ are semi-multiplicative.

\vspace{1mm}

Define the strongly continuous one-parameter
representation $\gamma$ on $A$ such that $\gamma_t = \al_t \, \be_t$ for every $t \in \R$.

\medskip

\begin{lemma}
We have the following properties :
\begin{itemize}
\item $\al^l$ and $\be^l$ commute
\item $\al^r$ and $\be^r$ commute
\item $\gam$ is semi-multiplicative and
$\gam^l_t = \al^l_t \, \be^l_t$  and $\gam^r_t = \al^r_t \, \be^r_t $ for every $t \in
\R$.
\end{itemize}
\end{lemma}

The proof of this lemma is straightforward and will therefore be left out.

\begin{proposition}
Consider $z \in \C$. Then $\overline{\al}_z \, \overline{\be}_z \subseteq
\overline{\gam}_z$, $\overline{\be}_z \, \overline{\al}_z \subseteq
\overline{\gam}_z$ and $D(\al_z \, \be_z) \cap D(\be_z \, \al_z)$ is a strict core
for $\overline{\gam}_z$.
\end{proposition}
\begin{demo} Take $x \in D(\overline{\al}_z \, \overline{\be}_z)$.

\medskip

Choose $e \in D(\al^r_z \, \be^r_z)$. Because $x$ belongs to $D(\overline{\be}_z)$ and $e$
belongs to $D(\be^r_z)$, we know by proposition \ref{pa2.prop3} that $x\,e$ belongs to
$D(\be_z)$ and $\be_z(x\,e)
= \be_z(x) \be^r_z(e)$.

Because $\be_z(x)$ belongs to $D(\overline{\al}_z)$ and $\be^r_z(e)$ belongs to
$D(\al^r_z)$, this implies that $\be_z(x\,e)$ belongs to $D(\al_z)$ and $\al_z(\be_z(x\,
e)) = \al_z(\be_z(x)) \, \al^r_z(\be^r_z(e)) = \al_z(\be_z(x)) \, \gam^r_z(e)$. So we see
that $x\,e$ belongs to $D(\gam_z)$ and $\gam_z(x\,e) = \al_z(\be_z(x)) \,
\gam^r_z(e)$.

\medskip

Because $D(\al^r_z \, \be^r_z)$ is a core for $\gam^r_z$, we get by corollary
\ref{pa2.cor1} the existence of a net $(e_k)_{k \in K}$ in $D(\al^r_z \, \be^r_z)$ such
that $(e_k)_{k
\in K}$ and $\bigl(\gam^r_z(e_k)\bigr)_{k \in K}$ converge strictly to 1.

\medskip

By the first part of the proof, we know for every $k \in K$ that $x \, e_k$ belongs to $
D(\gam_z)$ and $\gam_z(x \, e_k) = \al_z(\be_z(x)) \, \gam^r_z(e_k)$. So we get that $(x
\, e_k)_{k \in K}$ converges strictly to $x$ and that $\bigl( \gam_z(x \, e_k)
\bigr)_{k \in K}$ converges strictly to $\al_z(\be_z(x))$.

Hence, the strict closedness of $\overline{\gam}_z$ implies that $x$ belongs to
$D(\overline{\gamma}_z)$ and $\gamma_z(x) = \al_z(\be_z(x))$.

\medskip

So we have proven that $\overline{\al}_z \, \overline{\be}_z \subseteq
\overline{\gam}_z$. In a similar way, one proves that
$\overline{\be}_z \, \overline{\al}_z  \subseteq
\overline{\gam}_z$.
The last statement follows from the fact that $D(\al_z \, \be_z) \cap D(\be_z \, \al_z)$
is a core for $\gam_z$ and the fact that $D(\gam_z)$ is a strict core for
$\overline{\gam}_z$.
\end{demo}

\bigskip

\section{Tensor products of one-parameter representations}

We will prove in these section some basic properties about the tensor product of two
one-parameter representations.

\medskip

Consider two Banach spaces $E$,$F$. Let $\tau$ be a norm on $E \od F$ such that $\|\tau(x
\ot y)\| \leq \|x\| \, \|y\|$ for $x \in E$ and $y \in F$. We denote the completion of $E
\od F$ with respect to $\tau$ by $E \ot_\tau F$. The norm on $E \ot_\tau F$ is denoted by
$\|.\|_\tau$.

\medskip

If $S \in \cB(E)$, $T \in \cB(F)$ and $S \od T$ is continuous with respect to $\tau$, we
denote the unique continuous linear extension of $S \od T$ to $E \ot_\tau F$ by $S
\ot_\tau T$.

\bigskip

In order to define the tensor product of two one-parameter representations and $E$ and $F$
resp., they need to behave well with respect to the norm $\tau$ on $E \od F$. Therefore we
will introduce the following terminology.

\begin{terminology}
Consider a strongly continuous one-parameter representation $\al$ on $E$ and a strongly
continuous one-parameter representation $\be$ on $F$. Then we call the pair $\al$,$\be$
compatible with respect to $\tau$ if we have for every $t \in \R$ that $\al_t \od \be_t$
is continuous with respect to $\tau$ and has norm less or equal than 1.
\end{terminology}

It is then clear that every $\al_t \od \be_t$ is isometric with respect to $\tau$.

\bigskip\medskip

Let us fix for a while  the following objects :
\begin{itemize}
\item a strongly continuous one-parameter representation $\al$ on $E$
\item a strongly continuous one-parameter representation $\be$ on $F$
\end{itemize}
such that $\al$,$\be$ is compatible with $\tau$.
Then we can give the following definition

\medskip

\begin{definition}
We define the strongly continuous one-parameter representation $\al \ot_\tau \be$ on $E
\ot_\tau F$ such that we have for every $t \in \R$ that $(\al \ot_\tau \be)_t$ is the
unique continuous linear extension of $\al_t \od \tau_t$.
\end{definition}

\bigskip

Now we have the following  natural result concerning the analytic continuation of $\al_z
\ot_\tau \be_z$.

\begin{proposition} \label{pa6.prop1}
Consider $z \in \C$. Then $\al_z \od \be_z$ is closable and its closure is equal to $(\al
\ot_\tau \be)_z$.
\end{proposition}
\begin{demo}
Looking at definition \ref{pa1.def1}, it is not difficult to see that
$\al_z \od \be_z \subseteq  (\al \ot_\tau \be)_z$
This implies that $\al_z \od \be_z$ is closable and that its closure is a restriction of
$(\al \ot_\tau \be)_z$

Define $C$ as the set of elements in $E$ which are analytic with respect to $\al$ and
define $D$ as the set of elements in $F$ which are analytic with respect to $\be$. Then
the following holds :
\begin{enumerate}
\item The set $C$ is dense in $E$, $C \subseteq D(\al_z)$ and $\al_t(C) \subseteq C$
for $t \in \R$.
\item The set $D$ is dense in $F$, $D \subseteq D(\be_z)$ and $\be_t(D) \subseteq D$
for $t \in \R$.
\end{enumerate}
So we see that $C \od D$ is dense in $E \ot_\tau F$, $C \od D \subseteq D((\al \ot_\tau
\be)_z)$ and $(\al \ot_\tau \be)_t(C \od D) \subseteq C \od D$ for $t \in \R$. Then
corollary \ref{pa1.cor3} implies that $C \od D$ is a core for $(\al \ot_\tau \be)_z$.

This implies that $(\al \ot_\tau \be)_z$ is the closure of $\al_z \od \be_z$.
\end{demo}

\medskip

\begin{corollary}
Consider $z \in \C$, $C$ a core for $\al_z$ and $D$ a core for $\be_z$. Then $C \od D$ is
a core for $(\al \ot \be)_z$.
\end{corollary}

\bigskip\bigskip

For the rest of this section, we will concentrate on the case of \cst-algebras and
semi-multiplicative one-parameter representations. So consider two \cst-algebras $A$,$B$
and a \cst-norm $\tau$ on $A \od B$. We will also fix the following objects :
\begin{itemize}
\item a semi-multiplicative strongly continuous one-parameter representation $\al$ on $A$
\item a semi-multiplicative strongly continuous one-parameter representation $\be$ on $B$
\end{itemize}
such that $\al$,$\be$ is compatible with resepct to $\tau$.

\bigskip

\begin{result}
We have that the pairs $\al^l$,$\be^l$ and $\al^r$,$\be^r$ are both compatible with
respect to $\tau$.
\end{result}
\begin{demo}
Choose $t \in \R$. It is easy to check that $(\al^l_t \od \be^l_t)(x) (\al_t \od \be_t)(y)
= (\al_t \od \be_t)(x y)$ for $x,y \in A \od B$.

We know that $\al_t \od \be_t$ is  bijective and isometric (this last one because of  the
compatiblity of $\al$,$\be$). Choose $x \in A \od B$. Then
\begin{eqnarray*}
\|(\al^l_t \od \be^l_t)(x)\|
& = & \sup \, \{ \, \|(\al^l_t \od \be^l_t)(x) \, z\| \mid z \in A \od B
\text{ such that } \|z\| \leq 1  \, \} \\
& = & \sup \, \{ \, \|(\al^l_t \od \be^l_t)(x) \, (\al_t \od \be_t)(y)\|
\mid y \in A \od B \text{ such that } \|(\al_t \od \be_t)(y)\| \leq 1  \, \} \\
& = & \sup \, \{ \, \|(\al^l_t \od \be^l_t)(x) \, (\al_t \od \be_t)(y)\|  \mid
y \in A \od B \text{ such that } \|y\| \leq 1  \, \} \\
& = & \sup \, \{ \, \|(\al_t \od \be_t)(x y)\|  \mid y \in A \od B
\text{ such that } \|y\| \leq 1  \, \} \\
& = & \sup \, \{ \, \|x y\|  \mid y \in A \od B \text{ such that } \|y\| \leq 1
\, \} = \|x\|
\end{eqnarray*}
So we see that $\al^l_t \od \be^l_t$ is isometric.

The statement about $\al^r$,$\be^r$ is proven in the same way.
\end{demo}

So we have the following strongly continuous one-parameter representations on $A \ot_\tau
B$ :
$$\al^l \ot_\tau \be^l \hspace{3cm} \al \ot_\tau \be \hspace{3cm} \al^r \ot \be^r$$

\medskip

\begin{result}
We have that $\al \ot_\tau \be$ is semi-multiplicative and that $(\al \ot_\tau \be)^l =
\al^l \ot_\tau \be^l$ and that $(\al \ot_\tau \be)^r = \al^r \ot_\tau \be^r$
\end{result}
\begin{demo}
Consider $t \in \R$. Then it follows readily that $(\al^l \ot_\tau \be^l)_t(x) \, (\al
\ot_\tau \be)_t(y) = (\al \ot_\tau \be)_t(x y)$ for $x,y \in A \od B$. This equality will
then by continuity hold for $x,y \in A \ot_\tau B$.

The equation involving $\al^r \ot_\tau \be^r$ is treated in the same way.
\end{demo}

\bigskip

\begin{proposition}
Consider $z \in \C$. Then $\oal_z \od \overline{\be}_z$ is strictly closable and its
strict closure is equal to $\overline{\al \ot_\tau \be}_{\,z}$ \ .
\end{proposition}
\begin{demo}
Using definition \ref{pa2.def1}, it is easily seen that $\oal_z \od \overline{\be}_z
\subseteq  \overline{\al \ot_\tau \be}_{\,z}$ \ . This implies that $\oal_z \od
\overline{\be}_z$ is strictly closable and that its strict closure is a restriction of
$\overline{\al \ot_\tau \be}_{\,z}$ \ .

We know by proposition \ref{pa6.prop1} that $D(\al_z) \od D(\be_z)$ is a norm core for
$(\al \ot_\tau \be)_z$, so it is a strict core for $\overline{\al \ot_\tau \be}_{\,z}$ \ .
Hence the strict closure of $\oal_z \od \overline{\be}_z$ is equal to $\overline{\al
\ot_\tau \be}_{\,z}$ \ .
\end{demo}

\bigskip

\begin{remark} \rm
Under stronger conditions it is possible to prove generalizations of some of these
results. For instance is it possible to proof the following result. Consider Banach spaces
$E$,$F$ and a sub-cross norm $\tau$ on $E \od F$. Let $\al$ be a strongly continuous
one-parameter group on $E$ and $\be$ be a strongly continuous one-parameter group $F$ such
that we have for every $s,t \in \R$ that $\al_s \od \be_t$ is continuous with respect to
$\tau$ and has norm less or equal than 1.

Consider $y,z \in \C$. Then $\al_y \od \be_z$ is closable with respect to $\tau$ and you
could define $\al_y \ot_\tau \be_z$ to be the closure of $\al_y \od \al_z$.

\vspace{1mm}

A possible proof goes as follows. Define strongly continuous one-parameter representations
$\eta$ and $\th$ on $E \ot_\tau F$ such that $\eta_t = \al_t \ot_\tau \io$ and $\th_t =
\io \ot_\tau \be_t$ for $t \in \R$. Then $\eta$ and $\th$ commute. So $\eta_y \, \th_z$ is
closable by the previous section . We have also that $\al_y \od \be_z$ is a restriction of
$\eta_y \, \th_z$, so $\al_y \od \be_z$ is a closable.

\vspace{1mm}

All of this is of course related to two-parameter representations. You can define a
two-parameter representation $\gamma$ on $E \ot_\tau F$ such that $\gamma_{s,t} = \al_s
\ot \be_t$ for $s,t \in \R$. Then you can define in a similar way as for one-parameter
representations a closed mapping $\gamma_{y,z}$ in $E \ot_\tau F$. Then $\gamma_{y,z}$ is
the closure of $\al_y \od \be_z$.

\vspace{1mm}

Similar remarks apply of course to the material gathered in the previous section.
\end{remark}

\bigskip

\section{The Stone theorem for Hilbert C*-modules}

Consider a Hilbert C*-module $E$ over a C*-algebra $A$ and and a strongly continuous
unitary group  representation $u$ of $\R$ on $E$. Using the Stone theorem for strongly
continuous unitary group representations on a Hilbert space, we want to give a proof of
the existence of a unique strictly positive element $T \in \cR(E)$ such that $T^{is} =
u_s$ for every $s
\in \R$. Another proof of this result can be found in \cite{Ston}.

\medskip

Notice that $u$ is  a strongly continuous one parameter representation on $E$. So we can
use the theory of analytic continuations for $u$. As a consequence, we have for every $z
\in \C$ that $u_z$ is a  A-linear closed linear mapping from within $E$ into $E$ (the
$A$-linearity is easy to check, the rest are results from the theory of analytic
continuations).

\begin{lemma} \label{pa3.lem2}
Consider $z \in \R i$. Then $u_z \subseteq (u_z)^*$.
\end{lemma}
\begin{demo}
Choose $a,b \in D(u_z)$. Define the function $f,g$ from $S(z)$ into $E$ such that $f(y) =
u_y (a)$ and $g(y) = u_y (b)$ for every $y \in S(z)$. Then $f,g$ are analytic on $S(z)^0$
and continuous on $S(z)$. We have for every $t \in \R$ that
$$\lan f(t), b \ran - \lan a , g(-t) \ran = \lan u_t(a) , b \ran - \lan a , u_{-t}(b) \ran
= \lan u_t(a) , b \ran - \lan a , u_t^*(b) \ran =  \lan u_t(a) , b \ran  - \lan u_t(a) , b
\ran = 0 \ . $$
Because the function $S(z) \rightarrow  A : y \mapsto \lan f(y) , b \ran -
\lan a , g(-\overline{y}) \ran$
 is analytic on $S(z)^0$ and continuous on $S(z)$, this implies that
$$0 = \lan f(z) , b \ran - \lan a , g(-\overline{z}) \ran = \lan u_z(a) , b \ran  -
\lan a , u_{-\overline{z}}(b) \ran  \stackrel{(*)}{=} \lan u_z(a) , b \ran - \lan a ,
u_z(b) \ran \ ,$$
where in (*) we used the fact that $z \in \R i$.

So we get that $\lan u_z(a) , b \ran  = \lan a , u_z(b) \ran$. The lemma follows.
\end{demo}

\begin{lemma} \label{pa3.lem3}
Consider $z \in \R \, i$. Then the mapping $1+u_z$ has dense range.
\end{lemma}
\begin{demo}
The lemma is trivially true if $z=0$. Therefore suppose that $z \neq 0$.

Choose $\om \in E^*$ such that $\om=0$ on the range of $1+u_z$. Take an element $a \in E$
which is analytic with respect to $u$. Define the function $h$ from $\C$ into $\C$ such
that $h(y) = \om(u_y(a))$ for every $y \in \C$. Then $h$ is analytic on $\C$.

Choose $y \in \C$. Then $u_y(a)$ belongs to $D(u_z)$ and $u_z(u_y(a)) = u_{y+z}(a)$.
Because $\om=0$ on the range of $1+u_z$, this implies that
$$ \om(u_y(a)) + \om(u_{y+z}(a)) = \om\bigl(u_y(a) + u_z(u_y(a))\bigr)
= \om\bigl((1 + u_z)(u_y(a))\bigr) = 0 \ .$$
So we see that $h(y) = - h(y+z)$. \ \ \ \ (a)

From the theory of analytic continuations, we know that $h$ is bounded on each horizontal
strip. Therefore, (a) implies that $h$ is bounded on $\C$. By the Liouville theorem, we
get that $h$ is constant. Hence, using (a) once again, we see that $h(0) = - h(0+z) = -
h(0)$. Hence, $h(0) = 0$. Consequently, $\om(a) = 0$.

Because the analytic elements with respect to $u$ form a dense subset of $E$, it follows
that $\om=0$.

By the Hahn-Banach theorem, the lemma follows.
\end{demo}

\medskip

\begin{proposition}
Consider $z \in \R i$. Then $u_z$ belongs to $\cR(E)$.
\end{proposition}
\begin{demo}
From the theory of analytic continuations, we know that $u_z$ is a densely defined closed
$A$-linear mapping. By the first lemma of this section, we get that $(u_z)^*$ is also
densely defined.

Furhermore,
$$1+u_{2 z} =  1+ u_z \, u_z \subseteq 1 + (u_z)^* (u_z) \ ,$$
so the previuous lemma implies that $1 + (u_z)^* (u_z)$ has dense range.

The proposition follows from the definition on page 96 on \cite{Lan}.
\end{demo}

\medskip

\begin{corollary}
Consider $z \in \C$. Then $u_z$ belongs to $\cR(E)$.
\end{corollary}

This follows by the previous proposition because $u_z = \, u_{\text{\scriptsize Re}\,z}
\, u_{i \, \text{\scriptsize Im}\,z}$ and $u_{\text{\scriptsize Re}\,z}$ is unitary
(example 2 after lemma 2.4 of \cite{Wor6}).

\bigskip

The proof of the following result is copied from lemma 9.5 of \cite{Stra}.

\begin{proposition}
Consider $z \in \R i$. Then $u_z$ is a strictly positive element in $\cR(E)$.
\end{proposition}
\begin{demo}
Using lemma \ref{pa3.lem2} and proposition \ref{pa1.prop4}, we get that $u_z =
u_{\frac{z}{2}+\frac{z}{2}} = u_{\frac{z}{2}} \, u_{\frac{z}{2}} \subseteq
(u_{\frac{z}{2}})^* (u_{\frac{z}{2}})$. This implies that $\lan u_z \, v , v \ran \geq 0$
for $v \in D(u_z)$ \ \ \ (*)

The above inequality implies that $u_z$ is symmetric. We will now  prove that $u_z$ is
selfadjoint.

\vspace{1mm}

By (*), we have for every $v \in D(u_z)$ that
$$\lan (1+u_z) v , (1 + u_z) v \ran = \lan v ,v \ran
+ \lan u_z \, v ,v  \ran + \lan v , u_z v  \ran  + \lan u_z \, v , u_z \, v  \ran
\geq \lan v, v \ran $$
This implies that $1 + u_z$ is injective and that $(1 + u_z)^{-1}$ is a bounded linear
mapping (with norm $\leq 1$). As the inverse of the closed mapping $1 + u_z$, the mapping
$(1 + u_z)^{-1}$ is also closed. Combining these two facts, we get that $D\bigl((1 +
u_z)^{-1}\bigr)$ is closed. Hence lemma \ref{pa3.lem3} implies that $D\bigl((1 +
u_z)^{-1}\bigr) = E$.

Because $u_z$ is symmetric, $(1 + u_z)^{-1}$ is also symmetric.

\medskip

Choose $w \in D((u_z)^*)$. Take $x \in E$. Then $(1 + u_z)^{-1} x$ belongs to
$D(u_z)$. Hence
$$\lan (1 + u_z)^{-1}\bigl((u_z)^* w + w) , x \ran
= \lan (u_z)^* w + w , (1 + u_z)^{-1}  x \ran
= \lan w , (1 + u_z)\bigl((1 + u_z)^{-1} x \bigr) \ran
= \lan w , x \ran $$
So we see that $w = (1 + u_z)^{-1}\bigl((u_z)^* w + w)$, which implies that
$w \in D(u_z)$.

\medskip

Therefore we have proven that $u_z$ is positive. Corollary \ref{pa1.cor2} implies that
$u_z$ has dense range, so $u_z$ is strictly positive.
\end{demo}

\bigskip

\begin{remark}  \label{pa3.rem1} \rm
It is now time to get some help from Hilbert space theory. For this, we will use results
from section 11 of \cite{JK}. We will start with a small overview of the necessary
results.

\medskip

Consider $\om \in A^*_+$. Then we use the following notations and results in the sequel.
\begin{trivlist}
\item[\ \ $\bullet$] We have a Hilbert space $E_\om$ and a linear map $\la_\om$ from $E$
into $E_\om$ such that
\begin{enumerate}
\item $\la_\om(E)$ is dense in $E_\om$
\item We have that $\lan \la_\om(a) , \la_\om(b) \ran = \om(\lan a , b\ran)$ for every
$a,b \in E$.
\end{enumerate}
So $\la_\om$ is continuous and has norm $\leq 1$.

\medskip

The construction of this Hilbert space and this mapping is completely similar to the
GNS-construction of a positive functional on a \cst-algebra.

\medskip

\item[\ \ $\bullet$] Consider $S \in \cR(E)$. Then we have a densely defined closed
operator $S_\om$ in $E_\om$ such that $\la_\om(D(S))$ is a core for $S_\om$ and such that
$S_\om \la_\om(a) = \la_\om(S(a))$ for every $a \in D(S)$ (definition 11.5 of \cite{JK}).

It is then easy to see that $\la_\om(C)$ is a core for $S_\om$ if $C$ is a core for $S$.

\vspace{1mm}

If $S$ belongs to $\cL(E)$, we have  that $S_\om \in \cB(E_\om)$ and that $S_\om \,
\la_\om(a) = \la_\om(S(a))$ for every $a \in E$.

\medskip

The mapping $\cL(E) \rightarrow \cB(E_\om) : S \mapsto S_\om$ is a $^*$-homomorphism which
is strongly continuous on bounded subsets.

\medskip

\item[\ \ $\bullet$] We have for every $S \in \cR(E)$ that $(S_\om)^* = (S^*)_\om$.

\medskip

\item[\ \ $\bullet$] Consider a strictly positive element $S$ in $\cR(E)$. Then $S_\om$
is a strictly positive element in $\cB(E_\om)$ and $(S^z)_\om = (S_\om)^z$ for $z \in \C$
(result 11.8 and 11.14 of \cite{JK}).
\end{trivlist}

\medskip

\item[\ \ $\bullet$]

We have moreover the following separation property (proposition 11.16 of \cite{JK}).

Consider $R,S \in \cR(E)$. Then $R = S$ $\Leftrightarrow$ We have for every state $\om$ on
$A$ that $R_\om = S_\om$.
\end{remark}

\bigskip\medskip

Consider $\om \in A^*_+$. It is now clear that the mapping $\R \rightarrow \cB(E_\om) : s
\mapsto (u_s)_\om$  is a strongly continuous unitary group representation of $\R$ on
$E_\om$. By the Stone theorem for Hilbert spaces, there exists an injective positive
operator $\al_\om$ in $E_\om$ such that $(u_s)_\om = (\al_\om)^{is}$ for every $s \in \R$.

\begin{lemma} \label{pa3.lem1}
We have for every $z \in \C$ and every $\om \in A^*_+$ that $(u_z)_\om = (\al_\om)^{i z}$.
\end{lemma}
\begin{demo}
For every $n \in \N$ and $a \in E$, we define the element
$$a(n) = \frac{n}{\sqrt{\pi}} \int \exp(-n^2 t^2) \, u_t(a) \, dt \in E \ ,
$$
then $a(n)$ belongs to $D(u_z)$ and
$$u_z(a(n)) = \frac{n}{\sqrt{\pi}} \int \exp(-n^2 (t-z)^2) \, u_t(a) \, dt \ .$$
We also have that the set $\langle \, a(n) \mid a \in E, n \in \N \, \rangle$ is a core
for $u_z$ (see result \ref{pa1.res3}).

By remark \ref{pa3.rem1} , we see that the set $\langle \, \la_\om(a(n))  \mid a \in E, n
\in \N \, \rangle$ is a core for $(u_z)_\om$ \ \ \ (a) \newline  and
\begin{eqnarray*}
(u_z)_\om \, \la_\om(a(n))  & = & \la_\om\bigl(u_z(a(n))\bigr) = \frac{n}{\sqrt{\pi}} \int
\exp(-n^2(t-z)^2) \, \la_\om(u_t(a)) \, dt \\
& = & \frac{n}{\sqrt{\pi}} \int \exp(-n^2 (t-z)^2) \, (u_t)_\om \, \la_\om(a)  \, dt
\text{\ \ \ \ \ \ (b)}
\end{eqnarray*}
for all $a \in E$, $n \in \N$.

\medskip

We have for every $a \in E$, $n \in \N$ that
\begin{eqnarray*}
\la_\om(a(n)) & = & \frac{n}{\sqrt{\pi}} \int \exp(-n^2 t^2) \, \la_\om(u_t(a)) \, dt  \\
& = & \frac{n}{\sqrt{\pi}} \int \exp(-n^2 t^2) \, (u_t)_\om \, \la_\om(a)  \, dt \ \, \\
& = & \frac{n}{\sqrt{\pi}} \int \exp(-n^2 t^2) \, (\al_\om)^{it} \, \la_\om(a) \, dt \ \,
\end{eqnarray*}
Because $\la_\om(E)$ is dense in $E_\om$, this equation  and result \ref{pa1.res3} imply
that the set $\langle \, \la_\om(a(n)) \mid a \in E, n \in \N \, \rangle$ is a core for
$(\al_\om)^{iz}$.
\ \ \ (c)

The same equation implies for all $a \in E$ and  $n \in \N$ that
\begin{eqnarray*}
(\al_\om)^{iz} \, \la_\om(a(n)) & = &
\frac{n}{\sqrt{\pi}} \int \exp(-n^2(t-z)^2) \, (\al_\om)^{it} \, \la_\om(a)  \, dt \\
& = & \frac{n}{\sqrt{\pi}} \int \exp(-n^2(t-z)^2) \, (u_t)_\om \, \la_\om(a) \, dt
= (u_z)_\om \, \la_\om(a(n)) \ ,
\end{eqnarray*}
where in the last equality, equation (b) was used.

Combining this equality with results (a) and (c), we get that
$(u_z)_\om = (\al_\om)^{iz}$.
\end{demo}

\medskip

Now we define $T = u_{-i}$ , so $T$ is a strictly positive element in $\cR(E)$.

\begin{proposition}
Consider $z \in \C$. Then $u_z = T^{i z}$.
\end{proposition}
\begin{demo}
Using lemma \ref{pa3.lem1} twice in a row, we have for every state $\om$ on $A$ that
$$(u_z)_\om = (\al_\om)^{iz} = ( (u_{-i})_\om )^{i z} = (T_\om)^{iz} =
(T^{iz})_\om \ . $$ This implies that $u_z = T^{iz}$ by remark \ref{pa3.rem1}.
\end{demo}

\begin{proposition}
The element $T$ is the unique strictly positive element in $\cR(E)$ such that $u_t
= T^{it}$ for every $t \in \R$.
\end{proposition}
\begin{demo}
The existence has already been established. Let us prove quickly the uniqueness.
Therefore choose a strictly positive element $S$  in $\cR(E)$ such that $u_t = S^{it}$
for every $t \in \R$.

Take a state $\om$ on $A$. In this case, we have for every $t \in \R$ that
$$(S_\om)^{it} = (S^{it})_\om = (T^{it})_\om = (T_\om)^{it}$$ which implies that
$S_\om = T_\om$.

Hence $S = T$ by remark \ref{pa3.rem1}.
\end{demo}

\bigskip

So we have proven the following Stonian theorem in the previous part of this section.

\begin{theorem} \label{pa3.thm1}
Consider a Hilbert-\cst-module $E$ over a \cst-algebra $A$ and a strongly continuous
unitary group representation $u$ of $\R$ on $E$. Then there exists a unique strictly
positive element $T$  in $\cR(E)$ such that $u_t = T^{it}$. We have moreover that $u_z =
T^{iz}$ for every $z \in \C$.
\end{theorem}

\medskip

A immediate implication of this theorem is the following proposition.

\begin{proposition} \label{pa3.prop2}
Consider a Hilbert \cst-module $E$ over a \cst-algebra $A$ and a strictly positive element
$T$ in $\cR(E)$. Define the strongly continuous unitary group  representation of $\R$ on
$E$ such that $u_t = T^{it}$ for every $t \in \R$. Then $u_z = T^{iz}$ for every $t \in
\R$
\end{proposition}

\bigskip

We will also need the following consequence of theorem \ref{pa3.thm1}.

\medskip

First we need some terminology. According to section 10 of \cite{JK}, the set $\cR(B)$ can
be naturally embedded in $\cR(E)$. This embedding is denoted by the mapping $\cR(B)
\rightarrow \cR(E) : T \mapsto \tilde{T}$.

\vspace{1mm}

We will need the following properties of this embedding :
\begin{trivlist}
\item[\ 1.] We have for every $T \in M(B)$ that $\tilde{T} = T$
\item[\ 2.] Consider a strictly positive element $T \in \cR(B)$. Then $\tilde{T}$ is
also strictly positive and $\widetilde{T^z} = \tilde{T}^z$ for $z \in \C$.
\end{trivlist}

\medskip

\begin{proposition} \label{pa3.prop3}
Consider a Hilbert \cst-module $E$ over a \cst-algebra $A$ and a non-degenerate
sub-\cst-algebra $B$ of $\cL(E)$. Let $\al$ be an  strictly positive  element in $\cR(E)$.
Then there exists a strictly positive element $T$ affiliated with $A$ such that $\tilde{T}
= \al$
\vspace{2mm}

$\Leftrightarrow$ \begin{minipage}[t]{12cm}
\begin{enumerate}
\item We have for every $t \in \R$ and $b \in B$ that $\al^{it} \, b$ belongs to $B$.
\item The function $\R \rightarrow B : t \mapsto \al^{it} \, b$ is norm continous for
every $b \in B$.
\end{enumerate} \end{minipage}
\end{proposition}
\begin{demo}
If there exists a strictly positive element $T$ affiliated with $A$ such that $\tilde{T}
= \al$, we have for every $t \in \R$ that $\al^{it} = \tilde{T}^{it} = \widetilde{T^{it}}
= T^{it}$. So $\al$ satisfies the two statements above.

\medskip

Now suppose that the statements 1) and 2) hold. Define the mapping $u$ from $\R$ into
$\cL(E)$ such that $u_t = \al^{it}$ for every $t \in \R$. Statements 1) and 2) imply that
$u$ is a strongly continuous unitary group homomorphism from $\R$ into $M(B)$. By the
Stonian theorem, there exists a strictly positive element $T \, \eta \, B$ such that $u_t
= T^{it}$ for every $t \in \R$. This implies for every $t \in \R$ that $\tilde{T}^{it} =
\widetilde{T^{it}} = T^{it} = u_t = \al^{it}$, which implies that $\tilde{T} = \al$.
\end{demo}

\bigskip

\section{Implemented one-parameter representations}

Throughout this section, we will fix a Hilbert-\cst-module $E$ over a \cst-algebra $A$. We
will study one-parameter representations implemented by strictly positive elements in
$\cR(E)$.

We will in particular prove Hilbert-\cst-module versions of well-known Hilbert space
results (see proposition 9.24 of \cite{Stra}).

\bigskip

In this section, we will use the following notation (see section  2 of \cite{JK}).

\begin{notation}
Let $S$,$T$ be elements in $\cR(E)$ and $x$ an element in $\cL(E)$. Then we say that $x$
is a middle multiplier of $S$,$T$ $\Leftrightarrow$ We have that $x(\text{Ran}\,T)
\subseteq D(S)$, $S \, x \, T$ is bounded and $\overline{S \, x \, T}$ belongs to
$\cL(E)$.

If $x$ is a middle multiplier of $S$,$T$, we define $S \comp x \comp T = \overline{S \, x
\, T}$.
\end{notation}

\bigskip\bigskip

\begin{terminology}
Consider a non-degenerate sub-\cst-algebra $B$ of $\cL(E)$ and a strongly continuous
one-parameter representation $\al$ on $B$.
\begin{itemize}
\item Let $S$ and $T$ be two strictly positive
elements in $\cR(E)$. Then we say that $\al$ is implemented by $S$, $T$ if and only if
$\al_t(x) = S^{it} \, x \, T^{-it}$ for every $x \in B$ and $t \in \R$.
\item Let $S$ be a strictly positive element in $\cR(E)$. Then we say that $\al$ is
implemented by $S$ if and only if $\al_t(x) = S^{it}\, x \, S^{-it}$ for every $x \in B$
and $t \in \R$.
\end{itemize}
\end{terminology}

\medskip

\begin{proposition}
Consider a non-degenerate sub-\cst-algebra $B$ of $\cL(E)$ and a strongly continuous
one-parameter representation $\al$ on $B$ which is implemented by two strictly positive
elements $S$,$T$ in $\cR(E)$. Then $\al$ is semi-multiplicative,  $\al^l$ is implemented
by $S$ and $\al^r$ is implemented by $T$.
\end{proposition}
\begin{demo}
Fix for the moment $t \in \R$. We have for every $x,y \in B$ that
$$S^{it} \, (x \, y) \, S^{-it} = (S^{it} \, x \, T^{-it})\,(T^{it}\,y\,S^{-it}) =
(S^{it}\,x\,T^{-it})\,(S^{it}\,y^*\,T^{-it})^* = \al_t(x) \, \al_t(y^*)^* \ . $$ This
implies immediately that $S^{it}\,B\,S^{-it} \subseteq B$.

\medskip

So we can define a mapping $\be$ from $\R$ into $\cB(B)$ such that $\be_t(x) = S^{it} \, x
\, S^{-it}$ for every $t \in \R$ and $x \in B$.
We get immediately that $\be$ is a representation of $\R$ on $E$.

By the formula above, we have that $\be_t(x\,y) = \al_t(x) \, \al_t(y^*)^*$ for every $x,y
\in B$ and $t \in \R$. This implies that $\be$ is strongly continuous.

So we see that $\be$ is a strongly continuous one-parameter representation on $B$ which by
definition is implemented by $S$. It is also clear that $\be_t(x) \, \al_t(y) =
\be_t(x\,y)$ for every $x,y \in B$.

In a similar way, we get the existence of a strongly continuous one-parameter
representation $\gam$ on $B$ which is implemented by $T$ and such that $\al_t(x)
\, \gam_t(y) = \al_t(x\,y)$ for every $x,y \in B$.
\end{demo}

\medskip

\begin{result}
Consider a non-degenerate sub-\cst-algebra $B$ of $\cL(E)$ and a strongly continuous
\newline one-parameter representation $\al$ on $B$ which is implemented by two
strictly positive
elements $S$,$T$ in $\cR(E)$.
Then we have for every $x \in M(B)$ and $t \in \R$ that
$$\al^l_t(x) = S^{it} \, x \, S^{-it} \hspace{2cm}
\al_t(x) = S^{it} \, x \, T^{-it} \hspace{2cm}
\al^r_t(x) = T^{it} \, x \, T^{-it}  $$
\end{result}
\begin{demo}
We have a bounded net $(x_j)_{j \in J}$ in $B$ such that $(x_j)_{j \in J}$ converges
strictly to $x$. Then $(\al_t(x_j))_{j \in J}$ is a bounded net which converges strictly
to $\al_t(x)$.

Because $B$ is non-degenerate in $E$, we get that $(x_j)_{j \in J}$ converges strongly to
$x$, so $(S^{it} \, x_j \, T^{-i t})_{j \in J}$ converges strongly to $S^{i t} \, x \,
T^{-i t}$. Using the non-degeneracy of $B$ once again, we get also that $(\al_t(x_j))_{j
\in J}$ converges strongly to  $\al_t(x)$.

Because $\al_t(x_j) = S^{it} \, x_j \, T^{-it}$ for every $j \in J$, we must have that
$\al_t(x) = S^{i t} \, x \, T^{-i t}$.

\medskip

The equalities concerning $\al^l$ and $\al^r$ are proven in the same way.
\end{demo}

\bigskip\medskip

We will make use of the notations and results introduced in remark \ref{pa3.rem1}.

\medskip

\begin{lemma} \label{pa5.lem2}
Consider $\om \in A^*_+$, $S,T \in \cR(E)$ and $x \in \cL(E)$ such that $x$ is a middle
multiplier of $S$,$T$. Then $x_\om$ is a middle multiplier of $S_\om$,$T_\om$ and $S_\om
\comp  x_\om  \comp  T_\om = (S  \comp  x  \comp  T)_\om$.
\end{lemma}
\begin{demo}
Choose $v \in D(S^*)$ and $w \in D(T)$. Then $\la_\om(v) \in D(T_\om)$, $\la_\om(w) \in
D((S_\om)^*)$ and $T_\om \, \la_\om(v) = \la_\om(T v)$ and $(S_\om)^* \, \la_\om(w) =
\la_\om(S^* w)$. Hence
$$\lan x_\om \, T_\om(\la_\om(v)) , (S_\om)^* \, \la_\om(w) \ran
= \lan x_\om \, \la_\om(T v) , \la_\om(S^* w) \ran
= \lan \la_\om(x \, T(v) ) , \la_\om(S^* w) \ran
= \om(\lan x \, T(v) , S^*w \ran)$$
Because $x$ is a middle multiplier of $S$,$T$, we get by definition that
$x \, T(v) \in D(S)$ and that $S(x \, T(v)) = (S  \comp  x  \comp  T) \, v$.
This implies that
\begin{eqnarray*}
& & \lan x_\om \, T_\om(\la_\om(v)) , (S_\om)^* \la_\om(w) \ran
 =  \om(\lan S(x \, T(v)) , w \ran)
= \om(\lan (S \comp  x  \comp  T) \, v , w \ran ) \\
& & \spat =   \lan \la_\om((S  \comp  x  \comp  T) \, v) , \la_\om(w) \ran
= \lan (S  \comp  x  \comp  T)_\om \, \la_\om(v) , \la_\om(w) \ran
\end{eqnarray*}

Because $\la_\om(D(T))$ is a core for $T_\om$ and $\la_\om(D(S^*))$ is a core for
$(S_\om)^*$, the above equation implies for  every $p \in D(T_\om)$ and $q \in
D((S_\om)^*)$ that $\lan x_\om \, T_\om(p) , (S_\om)^*(q) \ran = \lan (S  \comp  x  \comp
T)_\om \, p , q \ran$.

This implies easily that $x_\om(\text{Ran}\,T_\om) \subseteq D(S_\om)$ and that
$S_\om \, x_\om \, T_\om \subseteq (S  \comp  x  \comp  T)_\om$.
Hence we get by definition that $x_\om$ is a middle multiplier of $S_\om, T_\om$ and
$S_\om  \comp  x_\om  \comp  T_\om = (S  \comp  x  \comp  T)_\om$.
\end{demo}

\medskip

\begin{lemma} \label{pa5.lem1}
Consider a subset $K$ of $A_+^*$ which is separating for $A$. Let $S$,$T$ be elements in
$\cR(E)$ and $x,y$  elements in $\cL(E)$ such that we have for every $\om \in K$ that
$x_\om$ is a middle multiplier of $S_\om$,$T_\om$ and that $S_\om  \comp  x_\om  \comp
T_\om = y_\om$. Then $x$ is a middle multiplier of $S$,$T$ and $S  \comp  x  \comp  T =
y$.
\end{lemma}
\begin{demo}
Take $v \in D(T)$. Choose $w \in D(S^*)$.

Choose also $\om \in K$. Then $\la_\om(v) \in D(T_\om)$, $\la_\om(w) \in D((S_\om)^*)$ and
$T_\om \, \la_\om(v) = \la_\om(T v)$ and $(S_\om)^* \, \la_\om(v) = \la_\om(S^* v)$.

Because $x_\om$ is a middle multiplier of $S_\om$,$T_\om$, we have by definition that
$x_\om \, T_\om(\la_\om(v)) \in D(S_\om)$ and that
$S_\om\bigl(x_\om \, T_\om(\la_\om(v))\bigr) = y_\om \, \la_\om(v) = \la_\om(y v)$.

This gives us that
\begin{eqnarray*}
& & \om(\lan x \, T(v) , S^*(w) \ran) = \lan \la_\om(x\, T(v)) , \la_\om(S^* w) \ran
= \lan x_\om \, T_\om(\la_\om(v)) , (S_\om)^* \, \la_\om(w) \ran \\
& & \spat = \lan S_\om\bigl(x_\om \, T_\om(\la_\om(v))\bigr) , \la_\om(w) \ran
= \lan \la_\om(y v) , w \ran = \om(\lan y v , w \ran)
\end{eqnarray*}
So we see that $\lan x \, T(v) , S^* w \ran = \lan y v , w \ran$.

This gives us that $x \, T(v) \in D(S)$ and that $S(x\,T(v)) = y v$.

Hence we get by definition that $x$ is a middle multiplier of $S$,$T$ and that $S  \comp
x  \comp  T = y$.
\end{demo}

\medskip

\begin{proposition}
Consider a non-degenerate sub-\cst-algebra $B$ of $\cL(E)$ and a strongly continuous
one-parameter representation $\al$ on $B$ which is implemented by two strictly positive
elements $S$,$T$ in $\cR(E)$. Let $z$ be a complex number and consider $x \in M(B)$. Then
we have the following properties.
\begin{enumerate}
\item We have that $x \in D(\oal_z)$ $\Leftrightarrow$
$x$ is a middle multiplier of $S^{i z}, T^{-i z}$ and $S^{i z} \comp  x \comp  T^{-i z}
\in M(B)$.
\item If $x$ belongs to $D(\oal_z)$, then $\al_z(x) = S^{i z}  \comp x \comp  T^{-i z}$.
\end{enumerate}
\end{proposition}
\begin{demo}
\begin{enumerate}
\item Suppose that $x$ is a middle multiplier of  $S^{i z}$,$T^{-i z}$ and that $S^{i z}
\comp  x \comp  T^{-i z}$ belongs to $M(B)$.

We define for every $\om \in A^*_+$ and $v,w \in E$ the element $\om_{v,w} \in B^*$ such
that $\om_{v,w}(y) = \om(\lan y v , w \ran)$ for $y \in B$. So $\om_{v,w}(y) = \om(\lan y
v , w \ran)$ for $y \in M(B)$.

Put $K = \{ \, \om_{v,w} \mid v,w \in E, \om \in A^*_+ \, \}$. We have then clearly that
$K$ separates $M(B)$. Because $\al$ is implemented by $S$,$T$ , we have also that $\th
\al_t \in K$ for every $\th \in K$ and $t \in \R$.

\medskip

Fix $\om \in A^*_+$ and $v,w \in E$. Then lemma \ref{pa5.lem2} implies that $x_\om$ is a
middle multiplier of $(S^{iz})_\om$,$(T^{-iz})_\om$ and that $(S^{iz})_\om  \comp  x_\om
\comp  (T^{-iz})_\om = (S^{iz} \comp  x  \comp  T^{-iz})_\om$. This implies that $x_\om$
is a middle multiplier of $(S_\om)^{iz}$,$(T_\om)^{-iz}$ and that $(S_\om)^{iz}  \comp
x_\om  \comp  (T_\om)^{-iz} = (S^{iz} \comp  x  \comp  T^{-iz})_\om$.

By proposition 9.24 of \cite{Stra}, there exists a function $f : S(z) \rightarrow \C$ 
such that $f$ is bounded and continuous on $S(z)$, analytic on $S(z)^0$ and such that
\begin{enumerate}
\item We have for every $t \in \R$ that $f(t) = \lan (S_\om)^{it} \, x_\om \,
(T_\om)^{-it} \la_\om(v) , \la_\om(w) \ran$.
\item $f(z) = \lan \bigl((S_\om)^{iz}  \comp  x_\om  \comp  (T_\om)^{-iz}\bigr)
\la_\om(v) , \la_\om(w) \ran$.
\end{enumerate}

This implies for every $t \in \R$ that
$$f(t) = \lan \la_\om(S^{it} \, x \, T^{-it} \, v) , \la_\om(w) \ran
= \lan \la_\om(\al_t(x) \, v) ,\la_\om(w) \ran
= \om(\lan \al_t(x) \, v , w \ran) = \om_{v,w}(\al_t(x))$$
We have moreover that
\begin{eqnarray*}
f(z) & = &\lan (S^{iz}  \comp  x  \comp  T^{-iz})_\om \, \la_\om(v) , \la_\om(w) \ran
= \lan \la_\om((S^{iz}  \comp  x  \comp  T^{-iz}) \, v) , \la_\om(w) \ran \\
& = & \om(\lan (S^{iz}  \comp  x  \comp  T^{-iz})\, v , w \ran) = \om_{v,w}(S^{iz}  \comp
x  \comp  T^{-iz})
\end{eqnarray*}

Because $S^{iz}  \comp  x  \comp  T^{-iz}$ belongs to $M(B)$, proposition \ref{pa2.prop8}
implies that $x \in D(\oal_z)$ and that $\oal_z(x) = S^{iz}  \comp  x  \comp  T^{-iz}$.
\item Suppose that $x \in D(\oal_z)$.

Take $\om \in A_+^*$.

The function $S(z) \rightarrow M(B) : u \mapsto \al_u(x)$ is bounded and strictly
continuous. This implies that the function $S(z) \rightarrow \cL(E) : u \mapsto \al_u(x)$
is bounded and strongly continuous. In turn, this gives us that the function $S(z)
\rightarrow \cB(E_\om) : u \mapsto (\al_u(x))_\om$ is bounded and strongly continuous.

We have also that the function $S(z)^0 \rightarrow \cL(E) : u \mapsto \al_u(x)$ is
analytic which implies that the function $S(z)^0 \rightarrow \cB(E_\om) : u \mapsto
(\al_u(x))_\om$ is analytic

Furthermore, $(\al_t(x))_\om = (S^{it} \, x \, T^{-it})_\om = (S_\om)^{it} \, x_\om \,
(T_\om)^{-it}$ for all $t \in \R$.

Therefore proposition 9.24 of \cite{Stra} implies that $x_\om$ is a middle multiplier of $(S_\om)^{iz}$,
$(T_\om)^{-iz}$ and that $(S_\om)^{iz}  \comp  x_\om   \comp  (T_\om)^{-iz} =
(\al_z(x))_\om$. Hence $x_\om$ is a middle multiplier of $(S^{iz})_\om$, $(T^{-iz})_\om$
and $(S^{iz})_\om  \comp  x_\om  \comp  (T^{-iz})_\om = (\al_z(x))_\om$.

\vspace{1mm}

So lemma \ref{pa5.lem1} implies tht $x$ is a middle multiplier of $S^{iz}$,$T^{-iz}$ and
that $S^{iz}  \comp  x  \comp  T^{-iz} = \al_z(x)$.
\end{enumerate}
\end{demo}

By proposition \ref{pa2.prop9}, this result implies the following one.

\begin{corollary}
Consider a non-degenerate sub-\cst-algebra $B$ of $\cL(E)$ and a strongly continuous
one-parameter representation $\al$ on $B$ which is implemented by two strictly positive
elements $S$, $T$ in $\cR(E)$. Let $z$ be a complex number and consider $x \in B$. Then $x
\in D(\al_z)$ $\Leftrightarrow$ $x$ is a middle multiplier of $S^{i z}$,$T^{-i z}$ and
$S^{i z} \comp  x \comp  T^{-i z} \in B$.
\end{corollary}

\bigskip

\section{Appendix : Strict extensions}

In this appendix, we formulate a condition under which a mapping between two \cst-algebras
has a unique extension to the multiplier algebra.

\begin{definition}
Consider two \cst-algebras $A$ and $B$ and a linear mapping $\rho$ from
$A$ into $M(B)$.

We call $\rho$ strict
$\Leftrightarrow$ \ \ \begin{minipage}[t]{10cm} \begin{trivlist}
\item[1.] $\rho$ is norm continuous.
\item[2.] $\rho$ is strictly continuous on bounded subsets of $A$.
\end{trivlist} \end{minipage}
\end{definition}

\medskip

\begin{proposition} \label{app.prop1}
Consider two \cst-algebras $A$ and $B$ and a linear mapping $\rho$ from $A$ into $M(B)$
which is strictly continuous on bounded subsets of $A$. Then there exists a unique linear
mapping $\th$ from $M(A)$ into $M(B)$ which extends $\rho$ and is strictly continuous on
bounded subsets of $M(A)$. We define $\overline{\rho} = \th$.

We use the notation $\rho(a) = \overline{\rho}(a)$ for every $a \in M(A)$.
\end{proposition}
\begin{demo}
Take an approximate unit $(e_k)_{k \in K}$ for $A$.

Let $x$ be an element in $M(A)$. Then we have the following properties:
\begin{enumerate}
\item It is clear that $(x e_k)_{k \in K}$ is a bounded net in $A$ which
converges strictly to $x$.

\item Choose a net $(a_i)_{i \in I}$ which is bounded and converges strictly to $x$.

Then $(a_i - a_j)_{(i,j) \in I \times I}$ is a bounded net in $A$ which converges strictly
to 0. This implies that the net $(\rho(a_i - a_j))_{(i,j) \in I \times I}$ converges
strictly to 0, which implies that $(\rho(a_i) - \rho(a_j))_{(i,j) \in I \times I}$
converges strictly to 0. Therefore, $(\rho(a_i))_{i \in I}$ is strictly Cauchy.

This implies that $(\rho(a_i))_{i \in I}$ is strictly convergent in $M(B)$.

\item Choose nets $(a_i)_{i \in I}$, $(b_j)_{j \in J}$ which  are bounded and
converge strictly to $x$. Let $c$, $d$ be elements in $M(B)$ such that $(\rho(a_i))_{i \in
I}$ converges strictly to  $c$, $(\rho(b_j))_{j \in J}$ converges strictly to $d$.

Now $(a_i - b_j)_{(i,j) \in I \times J}$ is a bounded net in $A$ which converges strictly
to 0. This implies that $(\rho(a_i - b_j))_{(i,j) \in I \times J}$ converges strictly to
0. Therefore we get that $(\rho(a_i)-\rho(b_j))_{(i,j) \in I \times J}$ converges strictly
to 0. But we have also that $(\rho(a_i)-\rho(b_j))_{(i,j) \in I \times J}$ converges
strictly to $c-d$, which implies that $c-d=0$.

Consequently, $c = d$
\end{enumerate}

\medskip

These 3 properties allow us to define a mapping $\th$ from $M(A)$ into $M(B)$ such that we
have for every $a \in M(A)$ and every bounded net $(a_i)_{i \in I}$ in $A$ which converges
strictly to $a$ that $(\rho(a_i))_{i \in I}$ converges strictly to $\th(a)$.

It is easy to check that $\th$ is a linear mapping which extends $\rho$.
We are now going to prove that $\th$ is strictly continuous on bounded subsets of $M(A)$.
Therefore choose a bounded net $(x_j)_{j \in J}$ in $M(A)$ and
$x \in M(A)$ such that $(x_j)_{j \in J}$ converges strictly to $x$.

Take a neighbourhood $V$ of $\th(x)$ for the strict topology on $M(B)$. Because the strict
topology is locally convex, there exists a strictly closed neighbourhood $U$ of $\th(x)$
for the strict topology on $M(B)$ such that $U \subseteq V$.

\begin{itemize}
\item The net $(x_j  e_k)_{(j,k) \in J \times K}$ is a bounded net in $A$ which
converges strictly to $x$. By the definition of $\th$, this implies that $(\rho(x_j
e_k))_{(j,k) \in J \times K}$ converges strictly to $\th(x)$. But this implies the
existence of $j_0 \in J$ and $k_0 \in K$ such that $\rho(x_j  e_k)$ belongs to $U$ for
every $j \in J$ and $k \in K$ with $j \geq j_0$ and $k \geq k_0$ \ \ \ (*).

\item Fix $j \in J$ with $j \geq j_0$. Then $(x_j e_k)_{k \in K}$ is a bounded
net in $A$ which converges strictly to $x_j$. By definition, we get that
$(\rho(x_j e_k))_{k \in K}$ converges strictly to $\th(x_j)$.

By (*), we know that $\rho(x_j e_k)$ belongs to $U$ for every $k \in K$ with $k \geq k_0$.
Therefore, the strict closedness of $U$ implies that $\th(x_j)$ belongs to $U$, so
$\th(x_j)$ belongs to $V$.
\end{itemize}
Hence, we find that $(\th(x_j))_{j \in J}$ converges strictly to $\th(x)$.

\medskip

Consequently, we have proven the existence. The uniqueness is trivial.
\end{demo}

It is rather straightforward to check that this extension behaves well with respect to
algebraic operations, e.g. if $\rho$ is multiplicative, then $\overline{\rho}$ will be
multiplicative.

\medskip

\begin{result}
Consider two \cst-algebras $A$ and $B$ and a strict linear mapping $\rho$ from $A$ into
$M(B)$ which is strictly continuous on bounded subsets of $A$. Then $\overline{\rho}$ is
bounded and $\|\overline{\rho}\| = \|\rho\|$.
\end{result}
\begin{demo}
Choose $x \in M(A)$. Take an approximate unit $(e_k)_{k \in K}$ for $A$. Then $(x e_k)_{k
\in K}$ is a bounded net in $A$ which converges strictly to $x$. Therefore, $(\rho(x
e_k))_{k \in K}$ converges strictly to $\overline{\rho}(x)$. We have for every $k \in K$
that $\|\rho(x e_k)\| \leq \|\rho\| \, \|x e_k\| \leq \|\rho\| \, \|x\|$ , so
$\|\overline{\rho}(x)\| \leq \|\rho\| \, \|x\|$.

This implies that $\overline{\rho}$ is bounded and $\|\overline{\rho}\| \leq \|\rho\|$.
Because $\overline{\rho}$ is an extension of $\rho$, we have that $\|\overline{\rho}\| =
\|\rho\|$.
\end{demo}

\medskip

Using these two results, we get the following one.

\begin{notation}
Consider \cst-algebras $A$, $\!B$, $\!C$. Let $\rho$ be a strict linear mapping from $A$
into $M(B)$ and $\th$ a strict linear mapping from $B$ into $M(C)$. We define $\th \rho =
\overline{\th} \circ \rho$. Then $\th \rho$ is a strict linear mapping from $A$ into
$M(C)$ with $\|\th \rho\| \leq \|\th\| \, \|\rho\|$ and such that $(\th \rho)(a) =
\th(\rho(a))$ for every $a \in M(A)$.
\end{notation}

\bigskip

Let $A$ be a \cst-algebra, $a$ an element in $M(A)$ and $\om \in A^*$. Then $\om a$ and $a
\om$ will denote the elements in $A^*$ such that $(\om a)(x) = \om(a x)$ and $(a \om)(x) =
\om(x a)$ for every $x \in A$.

\bigskip

Now we give some important examples of strict linear mappings:
\begin{trivlist}
\item[\ \ 1. ] Consider a \cst-algebra $A$ and $\om \in A^*$. Then there
exist $\th \in A^*$ and $a \in A$ such that $\om = \th a$ (see \cite{Tay}). This implies
easily that $\om$ is strict and $\om(x) = \th(a x)$ for every $x \in M(A)$.
\item[\ \ 2. ] Let $A$,$B$ be two \cst-algebras and $\pi$ a non-degenerate
$^*$-homomorphism from $A$ into $M(B)$ \ (non-degeneracy means that the set
$\pi(A) B$ is dense in $B$). Then $\pi$ is strict and $\pi(a) (\pi(b) c)
= \pi(a b) c$ for every $a \in M(A)$, $b \in A$ and $c \in B$ \ (this last
equality follows from the multiplicativity of $\overline{\pi}$).

Usually, the extension of non-degenerate $^*$-homomorphisms is defined by
the last formula.
\item[\ \ 3. ] Let $A$,$B$ be two \cst-algebras and $\rho$ a completely
positive mapping from $A$ into $M(B)$. We know that $\rho$ is bounded (see \cite{Tak}).
Let $(e_k)_{k \in K}$ be an approximate unit of $A$. By theorem 5.6  of \cite{Lan}, $\rho$
will be strict $\Leftrightarrow$ $(\rho(e_k))_{k \in K}$ is strictly convergent (in
\cite{Lan}, strictness of completely positive mappings is defined as the latter property).
\end{trivlist}

\bigskip

Consider \cst-algebras $A$,$B$,$C$,$D$ and consider the minimal tensor products $A \ot B$
and $C \ot D$.

Let $\rho : A \rightarrow C$ and $\eta : B \rightarrow D$ be linear mappings. We say that
$\rho$ and $\th$ are strictly tensorable if $\rho \od \th$ is continuous and if its unique
continuous linear extension $\rho \ot \th$ to $A \ot B$ is strict.

It is then easy to check that $(\rho \ot \th)(x \ot y)  = \rho(x) \ot \rho(y)$ for $x \in
M(A)$ and $y \in M(B)$.

\bigskip

From the third example of strict linear mappings, we get that any two strict completely
positive linear mappings are strictly tensorable. This implies that also a continuous
linear functional and a strict completely positive linear mapping are strictly tensorable
(this fact can also be proven more directly for a non-degenerate $^*$-homomorphism and a
continuous linear functional (see \cite{Hol1}).

\end{document}